\documentclass{amsart}
\usepackage{amscd}
\newfont{\cyr}{wncyr10}

\newcommand{\map}[1]{\;\xrightarrow{#1}\;}
\newcommand{\iso}{\cong}
\newcommand{\mil}{\lim\limits_\leftarrow}
\newcommand{\dlim}{\lim\limits_\rightarrow}
\newcommand{\Hom}{\mathrm{Hom}}
\newcommand{\Aut}{\mathrm{Aut}}  
\newcommand{\loc}{\mathrm{loc}} 
\newcommand{\Gal}{\mathrm{Gal}} 
\newcommand{\Tw}{\mathrm {Tw}}

\newcommand{\Q}{\mathbf Q}
\newcommand{\Z}{\mathbf Z} 
\newcommand{\RR}{{\mathrm{R}}}

\newcommand{\co}{\mathcal O}
\newcommand{\p}{\mathfrak p} 
\newcommand{\q}{\mathfrak q}
\newcommand{\fa}{\mathfrak a} 
\newcommand{\ff}{\mathfrak f}
 
\newcommand{\LL}{\mathcal L}

\newcommand{\Sel}{\mathrm {Sel}} 
\newcommand{\CS}{\mathcal{S}}
\newcommand{\CC}{\mathcal{C}}
\newcommand{\str}{\mathrm {str}} 
\newcommand{\rel}{\mathrm {rel}}
\newcommand{\HH}{\mathcal H} 
\newcommand{\gm}{\mathfrak m}

\newcommand{\f}{f} 
 
\newcommand{\tors}{\mathrm {tor}} 
\newcommand{\divs}{\mathrm {div}} 
\newcommand{\res}{\mathrm {res}}
\newcommand{\cor}{\mathrm {cor}}
\newcommand{\defect}{\mathfrak D}

\newcommand{\glob}{\mathrm {glob}}
\newcommand{\aug}{\mathcal {J}}
\newcommand{\reg}{\mathcal {R}}
\newcommand{\Der}{\mathrm {Der}}

\input xy \xyoption{all}
          
\newcommand{\bC}{{\mathbf C}}
\newcommand{\cL}{{\mathcal L}}
\newcommand{\rmN}{{\mathrm N}}
\newcommand{\fp}{{\mathfrak p}}
\newcommand{\fps}{{{\mathfrak p}^*}}
\newcommand{\bQ}{{\mathbf Q}}
\newcommand{\ov}[1]{{\bar{#1}}}
\newcommand{\vtheta}{\kappa}

\newcommand{\Zp}{\Z_p}
\newcommand{\Qp}{\Q_p}
\newcommand{\R}{\mathcal{R}}
\newcommand{\J}{\mathcal{J}}
\newcommand{\bfa}{\mathbf{a}}
\newcommand{\bfb}{\mathbf{b}}

\newcommand{\rank}{\mathrm{rank}}
\newcommand{\anti}{\mathrm{anti}}
\newcommand{\cycl}{\mathrm{cycl}}
\newcommand{\univ}{\mathrm{univ}}
\newcommand{\Res}{\mathrm{Res}}

\begin{document}
\title{Anticyclotomic Iwasawa theory of CM elliptic curves}
\date{Final version. January 11, 2006}
\author[A. Agboola \and B. Howard]
{Adebisi Agboola \and Benjamin Howard\\ \\With an appendix by
Karl Rubin} 

\address{Department of Mathematics\\ University of California\\ Santa
  Barbara, CA\\ 93106}
\email{agboola\char`\@math.ucsb.edu}
\address{Department of Mathematics\\ Harvard University\\ Cambridge, MA\\
02138}
\curraddr{Department of Mathematics\\ Boston College\\ Chestnut Hill, MA 02467}
\email{howardbe\char`\@bc.edu}
\address{Department of Mathematics, Stanford University, Stanford, CA
94305}
\curraddr{Department of Mathematics\\ UC Irvine\\ Irvine, CA 92697}
\email{krubin\char`\@math.uci.edu}

\thanks{The first author is partially supported by NSF
grant DMS-0070449}
\thanks{The second author is supported by a National Science Foundation 
Mathematical Sciences Postdoctoral Research Fellowship.}
\thanks{The author of the appendix is partially supported by NSF grant
DMS-0140378.}

\begin{abstract}
We study the Iwasawa theory of a CM elliptic curve $E$ in the
anticyclotomic $\mathbf{Z}_p$-extension of the CM field, where $p$ is
a prime of good, ordinary reduction for $E$.  When the complex
$L$-function of $E$ vanishes to even order, Rubin's proof of the two
variable main conjecture of Iwasawa theory implies that the Pontryagin
dual of the $p$-power Selmer group over the anticyclotomic extension
is a torsion Iwasawa module.  When the order of vanishing is odd, work
of Greenberg shows that it is not a torsion module.  In this paper we
show that in the case of odd order of vanishing the dual of the Selmer
group has rank exactly one, and we prove a form of the
Iwasawa main conjecture for the torsion submodule.
\end{abstract}

\maketitle

\theoremstyle{plain} 
\newtheorem{Thm}{Theorem}[subsection]
\newtheorem{Prop}[Thm]{Proposition} 
\newtheorem{Lem}[Thm]{Lemma}
\newtheorem{Cor}[Thm]{Corollary} 
\newtheorem{Conj}[Thm]{Conjecture}
\newtheorem{Theorem}{Theorem}
\newtheorem{bigConj}[Theorem]{Conjecture}

\theoremstyle{definition} 
\newtheorem{Def}[Thm]{Definition}

\theoremstyle{remark} 
\newtheorem{Rem}[Thm]{Remark}
\newtheorem{Ques}[Thm]{Question}
\setcounter{section}{-1}
\numberwithin{equation}{section}

\renewcommand{\theTheorem}{\Alph{Theorem}}

\section{Introduction and statement of results}

Let $K$ be an imaginary quadratic field of class number one, and let
$E_{/\Q}$ be an elliptic curve with complex multiplication by the
maximal order $\co_K$ of $K$. Let $\psi$ denote the $K$-valued
grossencharacter associated to $E$, and fix a rational prime $p>3$
at which $E$ has good, ordinary reduction.

Write $\Q_{p}^{\mathrm{unr}} \subset \mathbf{C}_p$ for the maximal
unramified extension of $\Q_p$, and let $R_0$ denote the completion of
its ring of integers. If $F/K$ is any Galois extension, then we write
$\Lambda(F)= \Z_p[[\Gal(F/K)]]$ for the generalised Iwasawa algebra,
and we set $\Lambda(F)_{R_0} = R_0[[\Gal(F/K)]]$. Let $C_{\infty}$ and
$D_{\infty}$ be the cyclotomic and anticyclotomic $\Z_p$-extensions of
$K$, respectively, and let $K_{\infty} = C_{\infty} D_{\infty}$ be the
unique $\Z_{p}^{2}$-extension of $K$.

 As $p$ is a prime of ordinary reduction for $E$, it follows that $p$
 splits into two distinct primes $p\co_K = \p \p^*$ over $K$. A
 construction of Katz gives a canonical measure
$$
\mathcal{L}\in\Lambda(K_\infty)_{\RR_0},
$$
the two-variable $p$-adic $L$-function, denoted 
$\mu_{\p^*}(K_\infty,\psi_{\p^*})$ in the text,
which interpolates the  value at $s=0$ of twists of $L(\psi^{-1},s)$ 
by characters of $\Gal(K_\infty/K)$.  It is a theorem of
Coates \cite{coates} that
the Pontryagin dual of the Selmer group
$\Sel_{\p^*}(E/K_\infty)\subset H^1(K_\infty,E[\p^{*\infty}])$ 
is a torsion $\Lambda(K_\infty)$-module, and
a fundamental theorem of Rubin, the two-variable 
Iwasawa main conjecture, asserts that
the characteristic ideal of this torsion module is generated by $\mathcal{L}$.
In many cases this allows one to deduce properties of the $\p^*$-power
Selmer group of $E$ over subfields of $K_\infty$.  For example,
if we identify 
$$
\Lambda(K_\infty)\iso\Lambda(D_\infty)
[[\Gal(C_\infty/K)]]
$$ 
and choose a topological generator
$\gamma\in\Gal(C_\infty/K)$, then we may expand $\mathcal{L}$
as a power series in $(\gamma-1)$
$$
\mathcal{L}=\mathcal{L}_{0}+
\mathcal{L}_{1}(\gamma-1)+ \mathcal{L}_{2}(\gamma-1)^2+\cdots 
$$
with $\mathcal{L}_i\in\Lambda(D_\infty)_{\RR_0}$.  
Standard ``control theorems'' imply that 
the characteristic ideal of 
$$
X^*(D_\infty)\stackrel{\mathrm{def}}{=}
\Hom\big(\Sel_{\p^*}(E/D_\infty),\Q_p/\Z_p\big)
$$  
(the reader should note Remark \ref{twisted action})
is then generated by the constant term $\mathcal{L}_0$.
If the sign in the functional equation of $L(E/\Q,s)$ is equal to $1$,
then theorems of Greenberg imply that $\mathcal{L}_0\not=0$, and
so $X^*(D_\infty)$ is a torsion 
$\Lambda(D_\infty)$-module.  In sharp contrast to this, 
when the sign of the functional equation 
is $-1$, the constant term vanishes and  $X^*(D_\infty)$
is not torsion, as was proved by Greenberg even before Rubin's proof
of the main conjecture (see \cite{greenberg}).

The following is the main result of this paper (which appears in the
text as Theorems \ref{main conjecture two} and \ref{T:bsd}) and was
inspired by conjectures of Mazur \cite{mazur}, Perrin-Riou \cite{BPR1}
and Mazur-Rubin \cite{MRtwo} concerning Heegner points.

\begin{Theorem}
Suppose the sign of the functional equation of $L(E/\Q,s)$ is $-1$.
Then $X^*(D_\infty)$
is a rank one $\Lambda(D_\infty)$-module.  If 
$\mathcal{X}\subset\Lambda(D_\infty)$ is the characteristic ideal
of the torsion submodule of $X^*(D_\infty)$, then
$$
\mathcal{X}\cdot\reg = (\mathcal{L}_1)
$$
as ideals of $\Lambda(D_\infty)_{R_0} \otimes_{\Z_p} \Q_p$, where 
$\reg$ is the regulator of the $\Lambda(D_\infty)$-adic height pairing
(defined in Section \ref{height section}).  
\end{Theorem}

A different statement of the Iwasawa main conjecture over $D_\infty$,
involving elliptic units and including the case where the sign
in the functional equation is equal to $1$, is also contained
in Theorem \ref{main conjecture two}.
Similar results in the Heegner point case
alluded to above can be found in \cite{me}.  T. Arnold \cite{arnold} has recently generalized Theorem A from the case of elliptic curves with complex multiplication to CM modular forms of higher weight.

The following result is due to K. Rubin. It establishes a conjecture
made in an earlier version of this paper, and a proof is given in the
Appendix. 

\begin{Theorem} \label{T:linterm}
Under the assumptions and notation of the Theorem A, the linear term
$\mathcal{L}_1$ is nonzero.
\end{Theorem}

While we have stated our results in terms of the $\p^*$-adic Selmer
group, they may equally well be stated in terms of the $\p$-adic
Selmer group.  If one replaces $\p^*$ by $\p$ in the above theorem,
then $\mathcal{X}$, $\mathcal{L}_1$, and $\reg$ are replaced by
$\mathcal{X}^\iota$, $\mathcal{L}_1^\iota$, and $\reg^\iota$,
respectively, where $\iota$ is the involution of $\Lambda(D_\infty)$
induced by inversion on $\Gal(D_\infty/K)$.  The decomposition
$E[p^\infty]\iso E[\p^\infty]\oplus E[\p^{*\infty}]$ of
$\Gal(K^\mathrm{al}/K)$-modules induces a decomposition of
$\Lambda(D_\infty)$-modules
$$
\Sel_{p}(E/D_\infty)\iso
\Sel_{\p}(E/D_\infty)\oplus \Sel_{\p^*}(E/D_\infty)
$$
which shows that, when the sign in the functional equation
is $-1$, the full $p$-power Selmer group 
$\Sel_{p}(E/D_\infty)$ has $\Lambda(D_\infty)$-corank
2, as was conjectured by Mazur \cite{mazur}.

An outline of this paper is as follows. The first section gives
definitions and fundamental properties of various Selmer groups
associated to $E$, with special attention to the anticyclotomic tower.
In the second section, we recall the definition of Katz's $p$-adic
$L$-function and the Euler system of elliptic units, and we state a
theorem of Yager which relates the two.  Our discussion of these
topics closely follows the excellent book of de Shalit
\cite{deShalit}.  Work of Rubin allows one to ``twist'' the elliptic
unit Euler system into an Euler system for the $\p$-adic Tate module
$T_\p(E)$, and we show, using nonvanishing results of Greenberg, that
the restriction of the resulting Euler system to the anticyclotomic
extension is nontrivial.  Applying the main results of \cite{rubin}
shows that a certain ``restricted'' Selmer group, contained in
$\Sel_{\p^*}(E/D_\infty)$, is a cotorsion module; using this we show
that $X^*(D_\infty)$ has rank one.  Using a form of Mazur's control
theorem, we then deduce that the characteristic ideal of a restricted
Selmer group over $K_\infty$ does not have an anticyclotomic zero.
This restricted two-variable Selmer group is related to elliptic units
by Rubin's proof of the two-variable main conjecture, and the
nonvanishing of its characteristic ideal along the anticyclotomic line
allows us to descend to the anticyclotomic extension and relate the
restricted Selmer group over $D_\infty$ to the elliptic units.  In the
third section we use results of Perrin-Riou and Rubin on the $p$-adic
height pairing to relate the twisted elliptic units to the linear term
$\mathcal{L}_1$ of Katz's $L$-function.
\smallskip

\noindent {\bf Acknowledgments.} We are very grateful to Karl Rubin
for writing an appendix to this paper.

\subsection{Notation and conventions}
We write $\psi$ for the $K$-valued grossencharacter associated to $E$,
and we let $\ff$ denote its conductor. Note that since $p$ is a prime
of good reduction for $E$, it follows that $p$ is coprime to $\ff$.

Let $\Q^\mathrm{al}\subset \mathbf{C}$ be the algebraic closure of
$\Q$ in $\mathbf{C}$, and let $\tau$ be complex conjugation, also
denoted by $z\mapsto\bar{z}$. Fix an embedding
$i_\p:\Q^\mathrm{al}\hookrightarrow\mathbf{C}_p$ lying above the prime
$\p$, and let $i_{\p^*}=i_\p\circ\tau$ be the conjugate embedding.

We write $R$ for the field of fractions of $R_0$. If $M$ is any
$\Z_p$-module, we define
$$
M_{\RR_0}=M\hat{\otimes}_{\Z_p}\RR_0,\hspace{1cm} M_\RR=M_{\RR_0}
\otimes_{\RR_0} \RR.
$$
The Pontryagin dual of $M$ is denoted
$$
M^\vee\stackrel{\mathrm{def}}{=}\Hom_{\Z_p}(M,\Q_p/\Z_p).
$$

If $M$ is any $\Z_p$-module of finite or cofinite type equipped with a
continuous action of $G_K=\Gal(\Q^\mathrm{al}/K)$ and $F$ is a
(possibly infinite) Galois extension of $K$, we let
$$
\HH^i(F,M)=\mil H^i(F',M)
$$
where the inverse limit is over all subfields
$F'\subset F$ finite over $K$ and is taken with respect to the natural
corestriction maps. If $\q$ is $\p$ or $\p^*$, we define
$$
\HH^i(F_\q,M)=\mil \bigoplus_{w|\q}H^i(F'_w,M),\hspace{1cm}
H^i(F_\q,M)=\dlim\bigoplus_{w|\q}H^i(F'_w,M).
$$ 
(Here the inverse (respectively, direct) limit is taken with
respect to the corestriction (respectively, restriction) maps.) These
groups have natural $\Lambda(F):= \Z_p[[\Gal(F/K)]]$-module
structures.

For a positive integer $n$, $C_n$ is the unique subfield of $C_\infty$
with $[C_n:K]=p^n$; the field $D_n$ is defined similarly.  If
$\mathfrak{m}$ is an ideal of $K$, we denote by $K(\mathfrak{m})$ the
ray class field of conductor $\mathfrak{m}$.  If $\mathfrak{n}$ is
another ideal of $K$, we let $K(\mathfrak{m}
\mathfrak{n}^\infty)=\cup_k K(\mathfrak{m}\mathfrak{n}^k)$. We write
$\mathrm{N}(\mathfrak{m})$ for the absolute norm of the ideal
$\mathfrak{m}$.


\section{A little cohomology}


We define canonical generators $\pi$ and $\pi^*$ of the ideals $\p$
and $\p^*$ by $\pi=\psi(\p)$ and $\pi^*=\psi(\p^*)$, so that
$\pi^*=\pi^\tau$.  Define $G_K$-modules
$$
W_\p=E[\p^\infty]\hspace{1cm} W_{\p^*}=E[\p^{*\infty}],
$$
and let $T_\p$ and $T_{\p^*}$ be the $\pi$ and $\pi^*$-adic Tate
modules, respectively.  Note that the action of $\tau$ on
$E[p^\infty]$ interchanges $W_\p$ and $W_{\p^*}$, and so induces a
group isomorphism $T_\p\iso T_{\p^*}$. If we set $V_\p=T_\p\otimes
\Q_p$ and $V_{\p^*}=T_{\p^*}\otimes \Q_p$, then there is an exact
sequence
\begin{equation}
\label{vector space sequence}
0\map{}T_\q\map{}V_\q\map{}W_\q\map{}0
\end{equation} 
where $\q=\p$ or $\p^*$. For every place $v$ of $K$ and any finite
extension $F$ of $K$ or $K_v$, the 
$G_F$-cohomology of this sequence implies that
$$
H^0(F,W_\p)_{/\divs}\iso H^1(F,T_\p)_\tors,
$$
where the subscript $/\divs$ indicates the quotient by the maximal
divisible submodule, and the subscript $\tors$ indicates the
$\Z_p$-torsion submodule.  The Weil pairing restricts to a perfect
pairing $T_\p\times T_{\p^*}\map{}\Z_p(1)$.


\subsection{Selmer modules}


Let $\q$ be either $\p$ or $\p^*$ and set $\q^*=\tau(\q)$.

\begin{Lem}\label{decomposition}
  The primes of $K$ above $p$ are finitely decomposed in $K_\infty$.
\end{Lem}
\begin{proof}
  This follows from Proposition II.1.9 of \cite{deShalit}.
\end{proof}

\begin{Lem}\label{no local torsion}
The degree of $K_\q(E[\q])$ over $K_\q$ is $p-1$.
\end{Lem}
\begin{proof}
This follows from the theory of Lubin-Tate groups.  See for example
Chapter 1 of \cite{deShalit}.
\end{proof}

\begin{Lem}\label{local finite}
For any intermediate field $K\subset F\subset K_\infty$ the $\Lambda(F)$-module
$$
A_\ff(F)=\bigoplus_{v|\ff}H^0(F_v,W_\p)
$$
has finite exponent.  If all primes dividing $\ff$ are finitely
decomposed in $F$, then $A_\ff(F)$ is finite.
\end{Lem}
\begin{proof}
Fix  a place $v|\ff$ of $F$. The extension $F_v/K_v$ is unramified, 
while $W_\p$ is
a ramified $\Gal(K_v^\mathrm{alg}/K_v)$-module (by the 
criterion of N\'eron-Ogg-Shafarevich).  Since $W_\p\iso\Q_p/\Z_p$
has no proper infinite submodules we conclude that $E(F_v)[\p^\infty]$
is finite, and so $\bigoplus_{v|\ff}H^0(F_v,W_{\p^*})$ has finite exponent.
This group is finite if all
primes above $\ff$ are finitely decomposed in $F$.
\end{proof}

\begin{Lem}\label{vector cohomology}
Let $F/K$ be a finite extension, let  $v$ be a prime of $F$ not
dividing $p$, and let $V$ be a finite  
dimensional $\Q_p$-vector space with a linear action of $G_{F_v}$.
If $V$ and $\Hom(V,\Q_p(1))$ both have no $G_{F_v}$-invariants then
$H^1(F_v,V)=0$.
\end{Lem}
\begin{proof}
This follows from Corollary 1.3.5 of \cite{rubin} and local duality,
or from standard properties of local Euler characteristics.
\end{proof}

In particular, Lemma \ref{vector cohomology} implies that
$H^1(F_v,V_\q)=0$ for $v$ not dividing $p$. If $v$ is any place of
$F$, we define the \emph{finite} or \emph{Bloch-Kato} local conditions
$$
H^1_\f(F_v,T_\q)=\left\{\begin{array}{ll}
    H^1(F_v,T_\q)_\tors& \mathrm{if\ }v\mid \q^*\\
    H^1(F_v,T_\q)&\mathrm{else}\end{array}\right.
$$
$$
H^1_\f(F_v,V_\q)=\left\{\begin{array}{ll}
    H^1(F_v,V_\q)\hspace{15pt} & \mathrm{if\ }v\mid \q\\
    0&\mathrm{else}\end{array}\right.
$$
$$
H^1_\f(F_v,W_\q)=\left\{\begin{array}{ll}
    H^1(F_v,W_\q)_\divs& \mathrm{if\ }v\mid \q\\
    0&\mathrm{else}\end{array}\right.
$$
The submodules
$H^1_\f(F_v,T_\q)$ and $H^1_\f(F_v,W_\q)$ are the preimage and image,
respectively, of $H^1_\f(F_v,V_\q)$ under the maps on cohomology
induced by the exact sequence (\ref{vector space sequence}).  If $M$
is any object for which we have defined $H^1_\f(F_v,M)$, we define the
\emph{relaxed} Selmer group $\Sel_\rel(F,M)$ to be the set of all
$c\in H^1(F,M)$ such that $\loc_v(c) \in H^1_\f(F_v,M)$ for every
place $v$ not dividing $p$.  We define the \emph{true} Selmer group
$\Sel(F,M)$ to be the subgroup consisting of all $c\in\Sel_\rel(F,M)$
such that $\loc_v(c)\in H^1_\f(F_v,M)$ for all $v$, including those
above $p$.  Finally, we define the \emph{strict} Selmer group to be
all those $c\in\Sel(F,M)$ such that $\loc_v(c)=0$ at the $v$ lying
above $p$.  By definition there are inclusions
$$
\Sel_\str(F,M)\subset\Sel(F,M)\subset\Sel_\rel(F,M)
$$
and all are $\Lambda(F)$-modules.  Our definitions of
$\Sel(F,T_\q)$ and $\Sel(F,W_\q)$ agree with the usual definitions of
the Selmer groups defined by the local images of the Kummer maps; see
Section 6.5 of \cite{rubin}.

\begin{Lem}\label{ramification}
Let $S$ denote the set of places of $K$ dividing $p\ff$ and let  $K_S/K$ be the maximal extension of $K$ unramified outside $S$. For any $K\subset F\subset K_S$ finite over $K$,
$$
H^1(K_S/F,T_\p)=\Sel_\rel(F,T_\p).
$$
\end{Lem}
\begin{proof}
For any $v\not\in S$, the local condition
$H^1_f(F_v, T_\p)$ is exactly the subgroup of unramified classes
by \cite{rubin} Lemma 1.3.5, while for $v\in S$ the local condition 
defining the relaxed Selmer group is all of $H^1(F_v,T_\p)$. 
\end{proof}

If $F/K$ is a (possibly infinite) extension we define
$\Lambda(F)$-modules
$$
\CS(F,T_\q)=\mil\Sel(F',T_\q) \hspace{1cm}
\Sel(F,W_\q)=\dlim\Sel(F',W_\q)
$$
where the limits are with respect to corestriction and restriction
respectively, and are taken over all subfields $F'\subset F$ finite
over $K$.  We also define strict and relaxed Selmer groups over $F$ in
the obvious way, e.g. $\CS_\str(F,T_\q)=\mil\Sel_\str(F',T_\q)$, and
so on.  Recall the notation
$$
\HH^i(F_\q,M)=\mil\bigoplus_{v|\q} H^i(F_v',M), 
$$
and set
$$
\HH_\f^1(F_\q,M)=\mil\bigoplus_{v|\q} H_\f^1(F_v',M), \quad
H_\f^1(F_\q,M)=\dlim\bigoplus_{v|\q} H_\f^1(F_v',M).
$$

The canonical involution of $\Lambda(F)$ which is inversion on
group-like elements is denoted $\iota:\Lambda(F)\map{}\Lambda(F)$.
This involution induces a functor from the category of
$\Lambda(F)$-modules to itself, which on objects is written as
$M\mapsto M^\iota$.  If $\q=\p$ or $\p^*$, there is a perfect local
Tate pairing
\begin{equation}\label{local pairing}
\HH^1(F_{\q},T_\p)\times H^1(F_{\q},W_{\p^*}) \map{}\Q_p/\Z_p
\end{equation}
which satisfies $(\lambda x,y)=(x,\lambda^\iota y)$ for
$\lambda\in\Lambda(F)$.  Under this pairing, the submodules
$\HH^1_\f(F_\q,T_\p)$ and $H^1_\f(F_\q,W_{\p^*})$ are exact orthogonal
complements.

\begin{Prop}\label{local ranks}
Suppose $F/K$ is a $\Z_p$-extension in which $\p^*$ ramifies.  Then
$\HH^1(F_{\p},T_\p)$ and $\HH^1(F_{\p^*},T_{\p})$ are rank one,
torsion-free $\Lambda(F)$-modules.
\end{Prop}
\begin{proof}
  Let $\q=\p$ or $\p^*$.  The claim that $\HH^1(F_\q,T_\p)$ has rank
  $1$ is Proposition 2.1.3 of \cite{perrin-riou}.  By Proposition
  2.1.6 of the same, the $\Lambda(F)$-torsion submodule of
  $\HH^1(F_\q,T_\p)$ is isomorphic to $H^0(F_\q,T_\p)$, so it suffices
  to show that $E(F_v)[\p^\infty]$ is finite for every place $v$ of
  $F$ above $\q$.  If $\q=\p$ this is immediate from Lemma \ref{no
  local torsion}.  If $\q=\p^*$ then $E[\p^\infty]$ generates an
  unramified extension of $K_\q$, and by hypothesis the intersection
  of this extension with $F_v$ is of finite degree over $K_\q$. Hence
  $E(F_v)[\p^\infty]$ is finite.
\end{proof}

\begin{Prop}\label{global duality}
If $F/K$ is an abelian extension such that 
$\HH^1_f(F_{\p^*},T_\p)=0$,  
then there are exact sequences
\begin{equation}\label{tate one}
0\map{}{\CS(F,T_\p)}\map{}{\CS_\rel(F,T_\p)}\map{\loc_{\p^*}}
\HH^1(F_{\p^*},T_\p)
\end{equation}
\begin{equation}
\label{tate two}
0\map{}{\Sel_\str(F,W_{\p^*})}\map{}{\Sel(F,W_{\p^*})}\map{\loc_{\p^*}}
H^1(F_{\p^*},W_{\p^*}),
\end{equation}
and the images of the rightmost arrows are exact orthogonal
complements under the local Tate pairing.  Under the same hypotheses
there are exact sequences
\begin{equation}
\label{tate three}
0\map{}{\CS_\str(F,T_\p)}\map{}{\CS(F,T_\p)}\map{\loc_\p}
\HH^1(F_{\p},T_\p)
\end{equation}
\begin{equation}
\label{tate four}
0\map{}{\Sel(F,W_{\p^*})}\map{}{\Sel_\rel(F,W_{\p^*})}\map{\loc_{\p}}
H^1(F_{\p},W_{\p^*}),
\end{equation}
and again the images of the rightmost arrows are exact orthogonal
complements under the sum of the local pairings. 
The hypotheses hold if $F$ contains a $\Z_p$-extension 
in which $\p^*$ ramifies.
\end{Prop}
\begin{proof}
If $\HH^1_f(F_{\p^*},T_\p)=0$ then 
  $H^1_\f(F_{\p^*},W_{\p^*})=H^1(F_{\p^*},W_{\p^*})$
by local duality.  The exactness
  of the sequences is now just a restatement of the definitions. The
  claims concerning orthogonal complements are consequences of
  Poitou-Tate global duality, c.f. Theorem 1.7.3 of \cite{rubin}.

Let $v$ be a place of $F$ above $\p^*$ and let 
$F'\subset F_v$ be finite over $K_{\p^*}$.  
From the cohomology of  (\ref{vector space sequence}) and the fact that
$H^0(F',V_\p)=0$ we have
$$
H^1_f(F',T_\p)\iso H^0(F',W_\p)\iso E(F')[\p^\infty].
$$
Taking the inverse limit over $F'\subset F_v$, we see that
$$\mil H^1_f(F',T_\p)=0$$ whenever $F_v$ contains an infinite pro-$p$
extension of $K_{\p^*}$ whose intersection with
$K_{\p^*}(E[\p^\infty])$ is of finite degree over $K_{\p^*}$.  The
extension of $K_{\p^*}$ generated by $E[\p^\infty]$ is unramified, so
this will be the case whenever $F$ contains a $\Z_p$-extension in which
$\p^*$ ramifies.
\end{proof}

If $F/K$ is an abelian extension we define $\Lambda(F)$-modules
\begin{eqnarray*}
X(F)&=&\Sel(F,W_\p)^\vee\\
X_\rel(F)&=&\Sel_\rel(F,W_\p)^\vee\\
X_\str(F)&=&\Sel_\str(F,W_\p)^\vee.
\end{eqnarray*}
Define $X^*(F)$, $X_\rel^*(F)$, and $X_\str^*(F)$ similarly,
replacing $\p$ by $\p^*$.

\begin{Rem}\label{twisted action}
Because of the behavior of the local pairing (\ref{local pairing})
under the action of $\Lambda(F)$, we adopt, for the entirety of the
paper, the convention that
$\Lambda(F)$ acts on $X(F)$ via $(\lambda\cdot f)(x)=f(\lambda^\iota x)$.
Thus the map
$$\HH^1(F_{\p^*},T_\p)\map{}X^*(F)$$
induced by localization at $\p^*$ and the local pairing is
a map of $\Lambda(D_\infty)$-modules.
The same convention is adopted for $X^*(F)$, $X_\rel(F)$, etc.
\end{Rem}

\begin{Lem}\label{equal ranks}
Let $F/K$ be a $\Z_p$ or $\Z_p^2$ extension of $K$.
There is a canonical isomorphism of $\Lambda(F)$-modules
$$
\CS_\rel(F,T_\p)\iso\Hom_{\Lambda(F)}
(X_\rel(F),\Lambda(F)).
$$
In particular 
$\CS_\rel(F,T_\p)$ and $X_\rel(F)$ have the same $\Lambda(F)$-rank,
and $\CS_\rel(F,T_\p)$ is torsion-free.
\end{Lem}
\begin{proof}
The proof is essentially the same as that of \cite{perrin-riou}
Proposition 4.2.3.
Suppose that $L\subset F$ is finite over $K$.
Let $S$ denote the set of places of $K$ consisting of the infinite place
and the prime divisors of $p\ff$, and let 
$K_S/K$ be the maximal extension of $K$ unramified outside 
$S$.  
By Lemma \ref{ramification}
$$
\CS_\rel(L,T_\p)\iso H^1(K_S/L,T_\p)\iso \mil H^1(K_S/L,E[\p^k]).
$$
On the other hand, Lemma \ref{no local torsion} shows that $E(L)[\p]=0$
(since $[L:K]$ is a power of $p$), and so
the $\Gal(K_S/L)$-cohomology of
$$
0\map{}E[\p^k]\map{}W_\p\map{\pi^k}W_\p\map{}0
$$
shows that $H^1(K_S/L,E[\p^k])\iso H^1(K_S/L,W_\p)[\p^k]$.

If we define
$$
X_S(L)=H^1(K_S/L,W_\p)^\vee
$$
then
$$\Hom_{\Lambda(L)}(X_S(L),\Lambda(L))\iso
\Hom_{\Z_p}(X_S(L),\Z_p)$$
via the augmentation map $\Lambda(L)\map{}\Z_p$. 
The right hand side is isomorphic to the $p$-adic Tate module of
$H^1(K_S/L,W_\p)$, so by the above
\begin{equation}\label{half way}
\CS_\rel(L,T_\p)\iso \Hom_{\Lambda(L)}(X_S(L),\Lambda(L)).
\end{equation}
As in the proof of Lemma \ref{ramification}, 
for any $v\not\in S$ the unramified classes in
$H^1(L_v,W_\p)$ agree with the local condition $H^1_f(L_v,W_\p)$,
and so we have the exact sequence
\begin{equation}\label{very relaxed}
0\map{}\Sel_\rel(L,W_\p)\map{} H^1(K_S/L,W_\p)\map{}
\bigoplus_{v|\ff}H^1(L_v,W_\p).
\end{equation}
The Pontryagin dual of the final term is
isomorphic to 
$$
\bigoplus_{v|\ff}H^1(L_v,T_{\p^*})\iso
\bigoplus_{v|\ff}H^0(L_v,W_{\p^*}),
$$
where we have used Lemma \ref{vector cohomology}
and the cohomology of the short exact sequence relating 
$T_{\p^*}$, $V_{\p^*}$, and $W_{\p^*}$.

If $A_\ff^*(F)$ denotes the module of Lemma \ref{local finite},
with $\p$ replaced by $\p^*$, we may take the limit as
$L$ varies, and the dual sequence to (\ref{very relaxed}) reads
$$ 
A_\ff^*(F)^\vee\map{}X_S(F)\map{}X_\rel(F)\map{}0
$$
where the first term is a torsion $\Lambda(F)$-module (even  a torsion 
$\Z_p$-module) by Lemma \ref{local finite}. Applying the functor
$\Hom_{\Lambda(F)}(\cdot,\Lambda(F))$ and
combining this with (\ref{half way})
gives the result.
\end{proof}

\begin{Rem}
A similar argument can be used to show that
$\CS(F,T_\p)$ and $X(F)$ have the same $\Lambda(F)$-rank, 
and similarly for the strict Selmer groups.  See \cite{perrin-riou}
Proposition 4.2.3, for example. 
When $F=D_\infty$, these facts will fall out 
during the more detailed analysis of the relationship between $X_\str$ and  
and $X_\rel$ given in Theorem \ref{strong selmer rank}.
\end{Rem}


\subsection{Anticyclotomic Iwasawa modules}
\label{anticyclotomic theory}


\begin{Lem}\label{change of group}
  There are isomorphisms of $\Lambda(D_\infty)$-modules
  $$
  \CS(D_\infty,T_\p)^\iota\iso\CS(D_\infty,T_{\p^*})\hspace{1cm}
  X(D_\infty)^\iota\iso X^*(D_\infty),
  $$
  and similarly for the relaxed and restricted Selmer groups.
\end{Lem}
\begin{proof}
  The action of $G_K$ on the full Selmer group
  $$\Sel(D_\infty, E[p^\infty])\iso \Sel(D_\infty,W_\p)\oplus
  \Sel(D_\infty,W_{\p^*})$$
  extends to an action of $G_\Q$, and
  complex conjugation interchanges the $\p$ and $\p^*$-primary
  components.  Since $\Gal(D_\infty/\Q)$ is of dihedral type, we may
  view complex conjugation as an isomorphism
  $$\Sel(D_\infty,W_\p)^\iota\iso \Sel(D_\infty,W_{\p^*}),$$
  and so
  $X(D_\infty)^\iota\iso X^*(D_\infty)$.  The other claims are proved
  similarly.
\end{proof}

The remainder of this subsection is devoted to a proof of the
following result.

\begin{Thm}\label{strong selmer rank}
If $r(\cdot)$ denotes $\Lambda(D_\infty)$-rank, then
$r(X(D_\infty))=r(\CS(D_\infty))$ and the same holds for 
the strict and relaxed Selmer groups.  Furthermore,
$$
r(X_\rel(D_\infty))=1+r(X_\str(D_\infty))
$$
and the $\Lambda(D_\infty)$-torsion submodules of  $X_\rel(D_\infty)$
and $X_\str(D_\infty)$ have the same characteristic ideals, up to powers
of $p\Lambda(D_\infty)$.
\end{Thm}

Let $\co$ be the ring of integers of some finite extension  $\Phi/\Q_p$,
and let 
$$
\chi:\Gal(D_\infty/K)\map{}\co^\times
$$ 
be a continuous character of $\Gal(D_\infty/K)$. If $M$ is any
$\Z_p$-module, define $M(\chi)=M\otimes\co(\chi)$. From Lemma \ref{no
local torsion} it follows that
$\psi_\q:G_K\map{}\Z_p^\times\map{}(\Z_p/p\Z_p)^\times$ is surjective,
and from this it is easy to see that the residual representation of
$T_\q(\chi)$ is nontrivial and absolutely irreducible.  Combining this
with Lemma \ref{vector cohomology} and the duality
$$
V_\p(\chi)\times V_{\p^*}(\chi^{-1})\map{}\Q_p(1),
$$
we have that $H^1(K_v,V_\q(\chi))=0$ for every $v$ not dividing
$p$.  We define generalized Selmer groups
$$
H^1_\rel(K,W_\q(\chi))\hspace{1cm}
H^1_\str(K,W_\q(\chi)),
$$ 
where the first group consists of classes which are everywhere
trivial at primes not dividing $p$, and which lie in the maximal
divisible subgroup of $H^1(K_v,W_\q(\chi))$ at primes $v$ above $p$,
and the second group consists of classes which are everywhere locally
trivial.

If $\chi$ is the trivial character,
then $H^1_\str(K,W_\q)=\Sel_\str(K,W_\q)$, but $H^1_\rel(K,W_\q)$
 may be slightly smaller than $\Sel_\rel(K,W_\q)$.
If we define 
$$\Sel_\rel(K,W_\q(\chi))\subset H^1(K,W_\q(\chi))$$
to be the subgroup of classes which are locally trivial at all primes
not dividing $p$, and impose no conditions at all above $p$
(so that this  agrees with our previous definition
when $\chi$ is trivial), then we can bound the index  of 
\begin{equation}\label{dumb relaxed condition}
H^1_\rel(K,W_\q(\chi))\subset \Sel_\rel(K,W_\q(\chi))
\end{equation}
as follows.  The quotient injects into
$$H^1(K_\p,W_\q(\chi))\oplus H^1(K_{\p^*},W_\q(\chi))$$
modulo its maximal divisible subgroup.  Thus, using 
the exact sequence (\ref{vector space sequence}) and local duality,
the order of the quotient is bounded by the order of
$$H^0(K_\p,W_{\q^*}(\chi^{-1}))
\oplus H^0(K_{\p^*},W_{\q^*}(\chi^{-1})).$$
It is easy to see that this group is finite, and bounded by some 
constant which does not depend on $\chi$, provided $\co$ remains fixed.

Our reason for working with the slightly smaller group $H^1_\rel$
is the following

\begin{Prop}(Mazur-Rubin)\label{MR}
For every character $\chi$ there is a non-canonical 
isomorphism of $\co$-modules
$$H^1_\rel(K,W_\p(\chi))\iso (\Phi/\co)\oplus 
H^1_\str(K,W_{\p^*}(\chi^{-1})).$$
\end{Prop}
\begin{proof}
All references in this proof are to \cite{mazur-rubin}.
It follows from Theorem 4.1.13 and Lemma 3.5.3 that
$$H^1_\rel(K,W_\p(\chi))[p^i]\iso (\Phi/\co)^r[p^i] \oplus 
H^1_\str(K,W_{\p^*}(\chi^{-1}))[p^i]$$
for every $i$, where $r$ is the core rank 
(Definition 4.1.11) of the local conditions
defining the Selmer group $H^1_\rel(K,W_\p(\chi))$.
A formula of Wiles, Proposition 2.3.5, shows that the core rank is
equal to 
$$\mathrm{corank}\ H^1(K_\p,W_{\p}(\chi))
+\mathrm{corank}\ H^1(K_{\p^*},W_{\p}(\chi))
-\mathrm{corank}\ H^0(K_v,W_{\p}(\chi)),
$$ 
in which $v$ denotes the unique archimedean place of $K$ and corank means corank as an $\co$-module.  The first two terms are each equal to $1$ by the local Euler characteristic formula, and the third is visibly $1$.  Hence the core rank is $1$.  Letting $i\to\infty$ proves the claim.
\end{proof}

Restriction gives a map
$$
H^1(K,W_\q(\chi))\map{}H^1(D_\infty,W_\q)(\chi)^{\Gal(D_\infty/K)},
$$
and since $H^0(D_\infty, W_\q(\chi))=0$, the Hochschild-Serre
sequence (see Proposition B.2.5 of \cite{rubin}) implies that
this map is an  isomorphism.

\begin{Lem}\label{twisted restriction}
The above restriction isomorphism induces injective maps
\begin{eqnarray*}
H^1_\str(K,W_\q(\chi))&\map{}&
\Sel_\str(D_\infty,W_\q)(\chi)^{\Gal(D_\infty/K)}\\
H^1_\rel(K,W_\q(\chi))&\map{}&
\Sel_\rel(D_\infty,W_\q)(\chi)^{\Gal(D_\infty/K)}
\end{eqnarray*}
whose cokernels are finite and bounded as $\chi$ varies (provided
$\co$ remains fixed).
\end{Lem}
\begin{proof}
A class $d\in \Sel_\str(D_\infty,W_\q)(\chi)^{\Gal(D_\infty/K)}$
is the restriction of some class $c\in H^1(K,W_\q(\chi))$
which is in the kernel of 
\begin{equation}\label{local restriction}
H^1(K_v,W_\q(\chi))\map{}H^1(D_{\infty,v},W_\q(\chi))
\end{equation}
for every place $v$ of $D_\infty$.
Let $\Gamma_v=\Gal(D_{\infty,v}/K_v)$, so that $\Gamma$ is either 
trivial or isomorphic to $\Z_p$.
If $\Gamma_v=0$ then (\ref{local restriction}) is an isomorphism.
If $\Gamma_v\iso\Z_p$ with generator $\gamma$, then the cokernel
is trivial and the kernel is isomorphic
to $M/(\gamma-1)M$ where $M=H^0(D_{\infty,v}W_\q)\otimes\co$.  
If $v$ is a prime
of good reduction not dividing $p$, then $D_{\infty,v}$ is the
unique unramified $\Z_p$-extension of $K_v$.  It follows that
 $M=0$ if $E[\q]\not\subset K_v$, while
$M=W_\q$ if $E[\q]\subset K_v$.  In either case,
since $(\gamma-1)$ acts as a nontrivial scalar on $W_\q$, we must
have $M/(\gamma-1)M=0$.  If $v$ is a prime of bad reduction, or if
$v$ lies above $p$, then $M$ is finite.
Thus the kernel of (\ref{local restriction}) is trivial for almost
all $v$, and finite and bounded by a constant independent of $\chi$.

A class $d\in \Sel_\rel(D_\infty,W_\q)(\chi)^{\Gal(D_\infty/K)}$
is the restriction of some class $c\in H^1(K,W_\q(\chi))$
which is in the kernel of (\ref{local restriction}) at every 
prime not dividing $p$.
The above argument shows that the cokernel of
$$
\Sel_\rel(K,W_\q(\chi))\map{}
\Sel_\rel(D_\infty,W_\q)(\chi)^{\Gal(D_\infty/K)}
$$ 
is finite with a bound of the desired sort, 
and so the claim follows from our bound on the index of 
(\ref{dumb relaxed condition}).
\end{proof}

\begin{Cor}
For any $\chi:\Gal(D_\infty/K)\map{}\co^\times$, the 
$\co$-coranks of 
$$\Sel_\rel(D_\infty,W_\p)(\chi)^{\Gal(D_\infty/K)}\hspace{1cm}
\Sel_\str(D_\infty,W_{\p})(\chi)^{\Gal(D_\infty/K)}$$
differ by $1$, and the quotients by the maximal $\co$-divisible 
submodules have the same order, up to $O(1)$ as $\chi$ varies.
\end{Cor}
\begin{proof}
By Lemma \ref{change of group},
$$
\Sel_\str(D_\infty,W_{\p^*})(\chi^{-1})^{\Gal(D_\infty/K)}
\iso
\Sel_\str(D_\infty,W_{\p})(\chi)^{\Gal(D_\infty/K)}.
$$
Combining this with Proposition \ref{MR} and 
Lemma \ref{twisted restriction} gives the stated result.
\end{proof}

\begin{Lem}\label{lambda core}
We have the equality $r(X_\rel(D_\infty))=1+r(X_\str(D_\infty))$,
and the $\Lambda(D_\infty)$-torsion submodules of  $X_\rel(D_\infty)$
and $X_\str(D_\infty)$ have the same characteristic ideals, up to powers
of $p\Lambda(D_\infty)$.
\end{Lem}
\begin{proof}
Choose a generator $\gamma\in\Gal(D_\infty/K)$ and
identify $\Lambda(D_\infty)$ with $\Z_p[[S]]$ via
$\gamma-1\mapsto S$.
Assume $\co$ is chosen large enough that the characteristic
ideals of the torsion submodules of $X_\str^*(D_\infty)$
and $X_\rel^*(D_\infty)$ split into linear factors.
Let $\gm\subset\co$ be the maximal ideal, and
fix pseudo-isomorphisms
$$X_\rel^*(D_\infty)\otimes_{\Z_p}\co\sim A\oplus A_p\hspace{1cm}
X_\str^*(D_\infty)\otimes_{\Z_p}\co\sim B\oplus B_p$$
where $A_p$ and $B_p$ are torsion modules with chracteristic ideals
generated by powers of $p$, and $A$ and $B$ are  of the form
$$A\iso\co[[S]]^a\oplus\bigoplus_{\xi\in\gm}A_\xi
\hspace{1cm}
B\iso\co[[S]]^b\oplus\bigoplus_{\xi\in\gm}B_\xi
$$
where each $A_\xi$ is isomorphic to $\bigoplus_i \co[[S]]/(S-\xi)^{e_{i}}$
for some exponents $e_i=e_i(A,\xi)$,
and similarly for $B$.  Define 
$$
P=\{\xi\in\gm\mid A_\xi \neq0\ \mathrm{or}\ B_\xi \neq 0\},
$$
and to any $\xi\not\in P$ we define a character $\chi_\xi$ by
$\chi_\xi(\gamma)=(\xi+1)^{-1}$.
Then for any $\xi\not\in P$ we have
$$
a=\mathrm{rank}_{\co}\ A/(S-\xi)A=\mathrm{corank}_\co\
\Sel_\rel(D_\infty,W_\p)(\chi_\xi)^{\Gal(D_\infty/K)}
$$
and similarly for $B$.
The corollary above immediately implies that $a=b+1$,
hence $r(X_\rel(D_\infty))=r(X_\str(D_\infty))+1$.

Now fix $\xi\in P$ and choose a sequence $x_k\to \xi$ with $x_k\in
{\mathfrak m}-P$ for all $k$.  As $k$ varies, the $\co$-length of the
torsion submodule of $A/(S-x_k)A$ is given by
$$ 
v(x_k-\xi)\cdot\sum_i e_i(A,\xi)\ +O(1)
$$ where $v$ is the 
valuation on $\co$, and similarly for $B$.  Applying the corollary,
we have 
$$
v(x_k-\xi)\cdot\sum_i e_i(A,\xi)=v(x_k-\xi)\cdot\sum_i e_i(B,\xi),
$$  
up to $O(1)$ as $k$ varies.  Letting $k\to\infty$ shows that 
$\sum_i e_i(A,\xi)=\sum_i e_i(B,\xi)$, proving that the torsion submodules 
of $X_\str(D_\infty)$ and $X_\rel(D_\infty)$ have the same
characteristic ideals, up to powers of $p\Lambda(D_\infty)$.
\end{proof}

The following corollary completes the proof of Theorem 
\ref{strong selmer rank}.

\begin{Cor}
We have the equality of ranks
$r(X(D_\infty))=r(\CS(D_\infty))$, and the same holds for 
the strict and relaxed Selmer groups.
\end{Cor}
\begin{proof}
For the relaxed Selmer groups, this equality of ranks was proved in
Lemma \ref{equal ranks}.  Let $A$ and $B$ be the cokernels of
$$
\CS_\rel(D_\infty,T_\p)\map{}\HH^1(D_{\infty,\p^*},T_\p),
  \hspace{1cm}\CS(D_\infty,T_\p)\map{}\HH^1(D_{\infty,\p},T_\p),
$$
respectively.  Then Propositions \ref{global duality} and \ref{change
of group} give
\begin{eqnarray*}
r(A)+r(X_\str(D_\infty))&=&r(X(D_\infty))\\
r(B)+r(X(D_\infty))&=&r(X_\rel(D_\infty))\\
r(A)+r(\CS_\rel(D_\infty,T_\p))&=&1+r(\CS(D_\infty,T_\p))\\
r(B)+r(\CS(D_\infty,T_\p))&=&1+r(\CS_\str(D_\infty,T_\p)).
\end{eqnarray*}
By Lemma \ref{lambda core}, the first two equalities imply that
$r(A)+r(B)=1$. The second two equalities then imply that
$r(\CS_\rel(D_\infty,T_\p))=1+r(\CS_\str(D_\infty,T_\p))$.
We deduce, using Lemma \ref{equal ranks}, that $r(X_\str(D_\infty))
=r(\CS_\str(D_\infty,T_\p))$.  Similarly, the equality 
$r(X(D_\infty))=r(\CS(D_\infty,T_\p))$ is deduced from 
Lemma \ref{equal ranks} by adding the middle two equalities.
\end{proof}

\section{$L$-functions and Euler systems}

In this section we recall the definition of Katz's $L$-function,
the construction of the elliptic units, and state Yager's theorem
relating the two.  Our presentation follows \cite{deShalit}, to which 
the reader is referred for more details on these topics.  Using results of
Rubin, we then twist the elliptic units into an 
Euler system more suitable for
our purposes and use the twisted Euler system to compute the corank of the $\p^*$-power Selmer group
over $D_\infty$.


\subsection{The $p$-adic $L$-function} 
\label{L functions}


For any integers $k,j$, we define a grossencharacter (of type $A_0$,
although we shall never consider any other type) of type $(k,j)$ to be
a $\Q^\mathrm{al}$-valued function, $\epsilon$, defined on integral
ideals prime to some ideal $\mathfrak{m}$, such that if
$\fa=\alpha\co_K$ with $\alpha\equiv 1\pmod{\mathfrak{m}}$ then
$\epsilon(\fa)=\alpha^k\bar{\alpha}^j$.  We have the usual notion of
the conductor of a grossencharacter, and the usual $L$-function
defined to be (the analytic continuation of)
$$
L(\epsilon,s)=\prod_\mathfrak{l}\frac{1}{
  1-\epsilon(\mathfrak{l})\mathrm{N} (\mathfrak{l})^{-s}}
$$
where the product is over all primes $\mathfrak{l}$ of $K$, with
the convention that $\epsilon(\mathfrak{l})=0$ for $\mathfrak{l}$
dividing the conductor of $\epsilon$.  For any ideal $\mathfrak{m}$,
the notation $L_{\mathfrak{m}}(\epsilon,s)$ means the $L$-function
without Euler factors at primes dividing $\mathfrak{m}$, and
$$
L_{\infty,\mathfrak{m}}(\epsilon,s)=
\frac{\Gamma(s-\min(k,j))}{(2\pi)^{s-\min(k,j)}}
L_\mathfrak{m}(\epsilon,s).
$$
Finally, if $\epsilon$ has conductor $\mathfrak{f}_\epsilon$ and type
$(k,j)$, set
$R(\epsilon,s)=(d_K\mathrm{N}(\mathfrak{f}_\epsilon))^{s/2}
L_{\infty,\mathfrak{f}_\epsilon}(\epsilon,s)$, where $d_K$ denotes the
discriminant of $K$. Then we have the functional equation
$R(\epsilon,s)=W_\epsilon\cdot R(\bar{\epsilon},1+k+j-s)$ for some
constant $W_\epsilon$ of absolute value one (the ``root number''
associated to $\epsilon$).  If we take $\epsilon=\psi$ to be the
grossencharacter of our elliptic curve, then the functional equation
reads $R(\psi,s)=W_\psi\cdot R(\psi,2-s)$, since
$\bar{\psi}(\mathfrak{l})=\psi(\bar{\mathfrak{l}})$ implies that
$L(\psi,s)=L(\bar{\psi},s)$.  In particular $W_\psi$ must be $\pm 1$.

To any grossencharacter of conductor dividing $\mathfrak{m}$, we
associate $p$-adic Galois characters
$$
\epsilon_\q:\Gal(K(\mathfrak{m}p^\infty)/K)\map{}\mathbf{C}_p^\times
$$
by the rule $\epsilon_\q(\sigma_\fa)=i_\q(\epsilon(\fa))$, where
$\q$ is $\p$ or $\p^*$, and $\sigma_\fa$ is the Frobenius of $\fa$.
The character $\psi_{\q}$ agrees with the character
$$\Gal(K(\ff \q^\infty)/K)\map{}\Aut(T_\q)\iso\Z_p^\times,$$
and the
formalism of the Weil pairing implies that $\psi_\p\psi_{\p^*}$ is the
cyclotomic character.

\begin{Thm}\label{Katz}(Katz)
  There are measures
\begin{equation} \label{E:interpolate 1}
  \mu_\p\in \Lambda(K(\ff p^\infty))_{\RR_0} \hspace{1cm}
  \mu_{\p^*}\in \Lambda(K(\ff p^\infty))_{\RR_0}
\end{equation}
  such that if $\epsilon$ is a grossencharacter of conductor
  dividing $\ff p^\infty$ of type $(k,j)$ with $0\le -j<k$, one has
  the interpolation formula
\begin{equation} \label{E:interpolate 2}
  \alpha_\p(\epsilon)\int \epsilon_\p\ d\mu_\p
  =\left(1-\frac{\epsilon(\p)}{p}\right) \cdot
  L_{\infty,\ff\p^*}(\epsilon^{-1},0)
\end{equation}
  where $\alpha_\p(\epsilon)\in\mathbf{C}_p$ is a nonzero constant,
  the integral is over $\Gal(K(\ff p^\infty)/K)$, and the right hand
  side is interpreted as an element of $\mathbf{C}_p$ via the
  embedding $i_\p$.  As usual, $\epsilon(\p)=0$ if $\p$ divides the
  conductor of $\epsilon$.  Similarly, if $\epsilon$ has infinity type
  $(k,j)$ with $0\le -j<k$ and conductor dividing $\ff p^\infty$, then
  $$
  \alpha_{\p^*}(\epsilon)\int\epsilon_{\p^*}\ d\mu_{\p^*}
  =\left(1-\frac{\epsilon(\p^*)}{p}\right) \cdot
  L_{\infty,\ff\p}(\epsilon^{-1},0)
  $$
  for some nonzero $\alpha_{\p^*}(\epsilon)$ where the right hand
  side is embedded in $\mathbf{C}_p$ via $i_{\p^*}$.
\end{Thm}
\begin{proof}
  This is Theorem II.4.14 of \cite{deShalit}.
\end{proof}

\begin{Rem}
  Our measure $\mu_\p$ is de Shalit's $\mu_p(\ff\p^{*\infty})$.  The
  measure $\mu_\p$ is canonically associated to the field $K$, the
  ideal $\ff$ and the embedding $i_\p$.  In particular it does not
  depend on the elliptic curve $E$. The constants
  $\alpha_\p(\epsilon)$ and $\alpha_{\p^*}(\epsilon)$ can be made
  explicit.
\end{Rem}

\begin{Rem} It can be deduced either from the interpolation formulae
  \eqref{E:interpolate 1}, \eqref{E:interpolate 2} or from a result of
   Yager (see Theorem \ref{yager} below) that the involution of
   $\Lambda(K(\ff p^\infty))_{\RR_0}$ induced by complex conjugation
   interchanges $\mu_\p$ and $\mu_{\p^*}$.
\end{Rem}

If $\epsilon$ is a grossencharacter of conductor dividing $\ff
p^\infty$, we define $$\LL_{\p,\ff}(\epsilon)= \int_{\Gal(K(\ff
  p^\infty)/K)}\epsilon^{-1}_\p\ d\mu_\p,$$
and similarly with $\p$
replaced by $\p^*$, so that the interpolation formula reads
\begin{equation}
\label{interpolation}
\LL_{\p,\ff}(\epsilon)=\left(1-\frac{\epsilon(\p)^{-1}}{p}\right)
\cdot L_{\infty,\ff\p^*}(\epsilon,0)
\end{equation}
up to a nonzero constant, provided that $\epsilon$ has infinity type
$(k,j)$ with $0\ge -j>k$.

If $\chi$ is a $\Z_p^\times$-valued character of $\Gal(F/K)$ for some
abelian extension $F/K$, we let
$$
\Tw_\chi:\Lambda(F)_{\RR_0}\map{} \Lambda(F)_{\RR_0}
$$
be the ring automorphism induced by
$\gamma\mapsto\chi(\gamma)\gamma$ on group-like elements.  Suppose
$\q=\p$ or $\p^*$, $F$ is an extension of $K$ contained in $K(\ff
p^\infty)$, and $\chi$ is a $\Z_p^\times$-valued character of
$\Gal(K(\ff p^\infty)/K)$.  Define $\mu_\q(F;\chi)$ to be the image of
$ \Tw_{\chi}(\mu_\q)$ under the natural projection
$$
\Lambda(K(\ff p^\infty))_{\RR_0}\map{} \Lambda(F)_{\RR_0},
$$
and for any integral $\co_K$-ideal $\fa$ prime to $\ff p$, define
$\lambda(F;\chi,\fa)$ to be the image of
$\Tw_{\chi}(\sigma_\fa-\mathrm{N}\fa)$. Let
$$
\mu_\q(F;\chi,\fa)= \mu_\q(F;\chi) \lambda(F;\chi,\fa).
$$
In particular, the measure $\mu_\q(D_\infty;\psi_{\p^*},\fa)$ will
be of crucial interest.

If $F/K$ is any subextension of $K(\ff p^\infty)$ and $\epsilon$ is a
grossencharacter such that $\epsilon_\q$ factors through $\Gal(F/K)$,
then
\begin{eqnarray}
\int\epsilon_\q^{-1}\ d\mu_\q(F;\psi_{\p^*},\fa)
&=&\big(\epsilon_\q^{-1}\psi_{\p^*}(\sigma_\fa)-\mathrm{N}\fa\big)
\int\epsilon_\q^{-1}\psi_{\p^*}\ d\mu_\q \nonumber \\
&=& \label{descent L value}
\left\{\begin{array}{ll}
 i_\p\big(\epsilon^{-1}\bar{\psi}(\fa)-\mathrm{N}\fa\big)\cdot
\LL_{\p,\ff}(\epsilon\bar{\psi}^{-1})& \mathrm{if\ }\q=\p\\ \\
i_{\p^*}\big(\epsilon^{-1}\psi(\fa)-\mathrm{N}\fa\big)\cdot
\LL_{\p^*,\ff}(\epsilon\psi^{-1})& \mathrm{if\ }\q=\p^*
\end{array}\right.
\end{eqnarray}
where the integral on the left is over $\Gal(F/K)$ and the integral on
the right is over $\Gal(K(\ff p^\infty)/K)$.

Fix a generator $\sigma$ of $\Gal(K_\infty/D_\infty)$.  The cyclotomic
character defines a canonical isomorphism
$$
< >:\Gal(K_\infty/D_\infty)\iso\Gal(C_\infty/K)\iso 1+p\Z_p.
$$
Following Greenberg, we define the \emph{critical divisor}
$$
\Theta=\sigma-<\sigma>\sigma^{-1}\in\Lambda(K_\infty)_{\RR_0}.
$$
If $\q$ is $\p$ or $\p^*$, we have a canonical factorization
$\psi_\q=\chi_\q\eta_q$ where $\chi_\q$ takes values in $\mu_{p-1}$
and $\eta_\q$ takes values in $1+p\Z_p$.  It is trivial to verify,
using the fact that $\eta_\p\eta_{\p^*}=< >$, that
$\Tw_{\eta_q}(\Theta)$ generates the kernel of the natural projection
$\Lambda(K_\infty)_{\RR_0}\map{} \Lambda(D_\infty)_{\RR_0}$.

Let $\mathcal{K}=K(E[p^\infty])$ so that $\mathcal{K}\subset K(\ff
p^\infty)$ by the theory of complex multiplication.  We have a natural
isomorphism
$$
\Gal(\mathcal{K}/K)\map{\psi_\p\times\psi_{\p^*}}
\Z_p^\times\times\Z_p^\times,
$$
and hence if we define $\Delta=\Gal(\mathcal{K}/K_\infty)$, every
character of $\Delta$ is of the form $\chi_\p^a\chi_{\p^*}^b$ for some
unique $0\le a,b<p-1$.  If $\chi$ is any character of $\Delta$, we let
$e(\chi)\in \Lambda(\mathcal{K})_{\RR_0}$ be the associated
idempotent, satisfying $\gamma e(\chi)=\chi(\gamma)e(\chi)$ We may
also view $e(\chi)$ as a map $\Lambda(\mathcal{K})_{\RR_0}\map{}
\Lambda(K_\infty)_{\RR_0}$, hopefully without confusion.  If $\chi$ is
the trivial character, this map is the natural projection.

\begin{Thm}\label{greenberg}(Greenberg)
  Let $\q=\p$ or $\p^*$, $\q^*=\bar{\q}$, and denote by $W$ the sign
  in the functional equation of $\psi$.  If $1$ denotes the trivial
  character, so that $\mu_\q(\mathcal{K},1)$ is the image of $\mu_\q$
  in $\Lambda(\mathcal{K})_{\RR_0}$, then
\begin{enumerate}
\item the critical divisor divides $e(\chi_\q)\mu_\q(\mathcal{K},1)$
  if and only if $W=-1$,
\item the critical divisor divides
  $e(\chi_{\q^*})\mu_\q(\mathcal{K},1)$ if and only if $W=1$.
\end{enumerate}
\end{Thm}
\begin{proof}
  The first claim is exactly the case $k_0=0$ of \cite{greenberg},
  Proposition 6.  For the second claim, the case $k_0=p-2$ of the same
  proposition shows that the critical divisor divides
  $e(\chi_{\q^*})\mu_\q(\mathcal{K},1)$ if and only if $W_{p-2}=-1$,
  where $W_{p-2}$ is the sign in the functional equation of
  $\psi^{2(p-2)+1}$. Let $m$ denote the number of roots of unity in
  $K$. Proposition 1 of \cite{greenberg}, together with the fact that
  $p\equiv 1\pmod{m}$, implies that $W_{p-2}=-W$.
\end{proof}

\begin{Cor}\label{parity}
  The measure $\mu_\p(D_\infty,\psi_{\p^*})$ is nonzero if and only if
  $W=-1$.  The measure $\mu_{\p^*}(D_\infty,\psi_{\p^*})$ is nonzero
  if and only if $W=1$.
\end{Cor}
\begin{proof} Let $\q=\p$ or $\p^*$.
  We have $\mu_\q(D_\infty,\psi_{\p^*})=0$ if and only if
  $\Tw_{\eta_{\p^*}}(\Theta)$ divides $\mu_\q(K_\infty,\psi_{\p^*})$,
  which occurs if and only if $\Theta$ divides
\begin{eqnarray*}
\Tw_{\eta_{\p^*}^{-1}}\big(\mu_\q(K_\infty,\psi_{\p^*})\big)
&=& \Tw_{\psi_{\p^*}^{-1}}\big(e(1)\mu_\q(\mathcal{K},\psi_{\p^*})\big)\\
&=&e(\chi_{\p^*})\Tw_{\psi_{\p^*}^{-1}}\big(
\mu_\q(\mathcal{K},\psi_{\p^*})\big)\\
&=&e(\chi_{\p^*})\mu_\q(\mathcal{K},1).
\end{eqnarray*}
The claim now follows from Theorem \ref{greenberg}.
\end{proof}

Corollary \ref{parity} shows that one of the measures
$\mu_{\fp}(D_{\infty}, \psi_{\fps})$, $\mu_{\fps}(D_{\infty},
\psi_{\fps})$ (depending upon the value of $W$) is non-zero. We
conclude this subsection by describing an alternative approach to
showing this fact, using root number calculations and non-vanishing
theorems for complex-valued $L$-functions.

Suppose that $\theta$ is a $\bC^{\times}$-valued idele class character
of $K$, of conductor $f_{\theta}$. Then we may write $\theta = \prod_v
\theta_v$, where the product is over all places of $K$. We define the
integer $n(\theta)$ by the equality $\theta_{\infty}(z) =
z^{n(\theta)} |z|^s$. Hence, if $\theta$ is associated to a
grossencharacter of $K$ of type $(k,j)$, then
$$
\theta_{\infty}(z) = z^k \ov{z}^j = z^{k-j} |z|^{2j},
$$
and so $n(\theta) = k-j$.

Let $W_{\theta}$ denote the root number associated to $\theta$. It
follows easily from standard properties of root numbers (see, for
example \cite{martinet}, especially Proposition
2.2 on page 30 and the definition on page 32) that $W_{\theta} =
W_{\theta / |\theta|}$.

\begin{Prop} (Weil) \label{P:weil}
Suppose that $\theta_1$ and $\theta_2$ are
$\bC^{\times}$-valued idele class characters of absolute value
$1$. Assume also that $f_{\theta_1}$ and $f_{\theta_2}$ are
relatively prime. Then
\begin{equation*}
W_{\theta_1} W_{\theta_2} \theta_1(f_{\theta_2}) \theta_2(f_{\theta_1})
=
\begin{cases}
W_{\theta_1 \theta_2} &\text{if $n(\theta_1) n(\theta_2) \geq 0$;} \\
(-1)^{\nu} W_{\theta_1 \theta_2}
&\text{if $n(\theta_1) n(\theta_2) < 0$,}
\end{cases}
\end{equation*}
where $\nu = \inf\{ |n(\theta_1)|, |n(\theta_1)|\} \mod{2}$.
\end{Prop}

\begin{proof} See \cite{weil}, pages 151-161 (especially Section 79).
\end{proof}

Let $K[p^n]$ denote the ring class field of $K$ of conductor $p^n$,
and set $K[p^{\infty}] = \cup_{n \geq 1} K[p^n]$.

\begin{Prop} \label{P:rootcalc}
Let $\xi$ be a grossencharacter of $K$ whose associated
Galois character factors through $\Gal(K[p^{\infty}]/K)$ and is of
finite order.

(a)(Greenberg) The following equality holds
$$
W_{\xi \psi} = W_{\psi}.
$$

(b) Let $e$ be a positive integer, and let $\epsilon$ be a
grossencharacter of type $(-e,e)$ associated to $K$, of trivial
conductor (such a grossencharacter always exists for a suitable choice
of $e$). Then
$$
W_{\epsilon \xi \psi} = -W_{\psi}.
$$
\end{Prop}

\begin{proof}
(a) This is proved on page 247 of \cite{greenberg}. (Note that the
proof given in \cite{greenberg} assumes that the Galois character
associated to $\xi$ factors through $\Gal(D_{\infty}/K)$. It is easy
to see that the same proof holds if we instead assume that the Galois
character associated to $\xi$ factors through $\Gal(K[p^{\infty}]/K)$.) 

(b) From part (a), we see that it suffices to show that
$$
W_{\epsilon \xi \psi} = -W_{\xi \psi}.
$$ 
The proof of this equality proceeds by applying Proposition
\ref{P:weil} to the idele class characters $\theta_1 = \xi \psi/|\xi
\psi|$ and $\theta_2 = \epsilon /|\epsilon|$.

We first note that since $\epsilon$ has trivial conductor, the same is
true of $\theta_2$, and so $\theta_1(f_{\theta_2}) =1$.  Next, we
observe that since $E$ is defined over $\bQ$, it follows that
$\ov{f}_{\theta_1} = f_{\theta_1}$. Since $\epsilon$ is of type
$(e,-e)$ and has trivial conductor, this implies that
$\theta_2(f_{\theta_1}) =1$.

Let $\delta_K$ denote the different of $K/\bQ$. It follows from
standard formulae for global root numbers (see \cite{lang}, Chapter
XIV, \S8, Corollary 1, for example) that
\begin{equation} \label{E:rootform}
W_{\epsilon} = i^{-2e} \epsilon(\delta_{K}^{-1}) = (-1)^e
\epsilon(\delta_{K}^{-1}).
\end{equation}
It is not hard to check that the different of any imaginary quadratic
field of class number one has a generator $\delta$ satisfying $\delta
/|\delta| = i$. This implies (since $\epsilon$ has trivial conductor
and is of type $(e,-e)$) that $\epsilon(\delta_{K}^{-1}) = (-1)^e$. It
now follows from \eqref{E:rootform} that
$$
W_{\epsilon} = 1 = W_{\theta_2}.
$$

Finally, we note that $n(\theta_1) = 1$ and $n(\theta_2) = -2e$,
whence it follows that $\nu = 1$. Putting all of the above together
gives
\begin{equation*}
W_{\theta_1 \theta_2} = -W_{\theta_1} W_{\theta_2},
\end{equation*}
from which we deduce that
\begin{equation*}
W_{\epsilon \xi \psi} = - W_{\xi \psi} =- W_{\psi} ,
\end{equation*}
as claimed.
\end{proof}

Now suppose that $W_{\psi} = -1$, and let $\vtheta = \epsilon \xi$ be a
grossencharacter as in Proposition \ref{P:rootcalc}(b) whose
associated Galois character factors through
$\Gal(D_{\infty}/K)$. (Note that once a choice of $\epsilon$ is fixed,
then there are infinitely many choices of $\xi$ such that the Galois
character associated to $\epsilon \xi$ factors through
$\Gal(D_{\infty}/K)$.) Then $\vtheta \ov{\psi}^{-1}$ is of type $(-e,
e-1)$, and so it lies within the range of interpolation of
\eqref{interpolation}. Hence we have (from \eqref{interpolation} and
\eqref{descent L value})
\begin{align*}
\int \vtheta_{\fp}^{-1} d \mu_{\fp}(D_{\infty};\psi_{\fps}, \fa) &=
i_{\fp}(\vtheta^{-1} \ov{\psi}(\fa) - \rmN \fa) \cL_{\fp,\ff}(\vtheta
\ov{\psi}^{-1}) \\
&= i_{\fp}(\vtheta^{-1} \ov{\psi}(\fa) - \rmN \fa)
\left(1 - \frac{\vtheta(\fp)^{-1}}{p} \right) L_{\infty, \ff
\fps}(\vtheta \ov{\psi}^{-1},0) \\
&= i_{\fp}(\vtheta^{-1} \ov{\psi}(\fa) - \rmN \fa)
\left(1 - \frac{\vtheta(\fp)^{-1}}{p} \right)
L_{\infty, \ff
\fps}(\vtheta \psi,1),
\end{align*}
where for the last equality we have used the fact that $\psi \ov{\psi}
= \rmN$.

Next we note that Proposition \ref{P:rootcalc}(b) implies that
$W_{\vtheta \psi} = 1$. It now follows from a theorem of Rohrlich
(see page 384 of \cite{rohrlich}) that, for all but finitely many choices of
$\vtheta$, we have $L_{\infty, \ff \fps}(\vtheta \psi,1) \neq
0$. Hence the measure $\mu_{\fp}(D_{\infty};\psi_{\fps}, \fa)$ is
non-zero, and so the same is true of
$\mu_{\fp}(D_{\infty};\psi_{\fps})$.

We now turn to the measure
$\mu_{\fps}(D_{\infty};\psi_{\fps})$. Suppose that $W_{\psi}=1$, and
let $\xi$ be any character of $\Gal(D_{\infty}/K)$ of finite
order. Then the grossencharacter $\xi \psi^{-1}$ is of type $(-1,0)$,
and so lies within the range of interpolation of
\eqref{interpolation}. Hence, just as above, we have
\begin{align*}
\int \xi_{\fps}^{-1} d \mu_{\fps}(D_{\infty};\psi_{\fps}, \fa) &=
i_{\fps}(\xi^{-1} \psi(\fa) - \rmN \fa) \cL_{\fps,\ff}(\xi
\psi^{-1}) \\
&= i_{\fps}(\xi^{-1} \psi(\fa) - \rmN \fa)
\left(1 - \frac{\xi(\fps)^{-1}}{p} \right) L_{\infty, \ff
\fp}(\xi \psi^{-1},0) \\
&= i_{\fps}(\xi^{-1} \psi(\fa) - \rmN \fa)
\left(1 - \frac{\xi(\fps)^{-1}}{p} \right)
L_{\infty, \ff
\fp}(\xi \ov{\psi},1).
\end{align*}

Now Proposition \ref{P:rootcalc}(a) implies that that $W(\xi \ov{\psi})
= W(\psi) = 1$ . It therefore follows from the theorem of Rohrlich
quoted above that, for all but finitely many choices of $\xi$, we have
$L_{\infty, \ff \fp}(\xi \ov{\psi},1) \neq 0$. This in turn implies
that $\mu_{\fps}(D_{\infty};\psi_{\fps}, \fa)$ is non-zero, whence it
follows that $\mu_{\fps}(D_{\infty};\psi_{\fps})$ is non-zero also.


\subsection{Elliptic units}


If $F/K$ is any finite extension and $\q=\p$ or $\p^*$, we define
$U_\q(F)$ to be the direct sum over all places $w$ dividing $\q$ of
the principal units of $F_w$.  If $F/K$ is any (possibly infinite)
extension, we define $U_\q(F)$ to be the inverse limit with respect to
the norm maps of the groups $U_\q(F')$ as $F'\subset F$ ranges over
the finite extensions of $K$.  Let $\fa$ be an $\co_K$-ideal with $(\fa,\ff
p)=1$, and let $I(\fa)$ denote the set of all ideals of $\co_K$ prime
to $\fa$.

If $L\subset\mathbf{C}$ is a lattice with CM by $\co_K$, we define an
elliptic function
$$\Theta(z;L,\fa)=\frac{\Delta(L)}{\Delta(\fa^{-1}L)}
\prod\frac{\Delta(L)}{(\wp(z,L)-\wp(u,L))^6}$$
where the product is
over the nontrivial $u\in\fa^{-1}L/L$, and $\Delta$ is the modular
discriminant.  For any $\mathfrak{m}\in I(\fa)$, let $w_\mathfrak{m}$
be the number of roots of unity congruent to $1$ modulo
$\mathfrak{m}$.  If $\mathfrak{m}$ is not a prime power then
$\Theta(1;\mathfrak{m},\fa)$ is a unit of $K(\mathfrak{m})$, and if
$\mathfrak{l}\in I(\fa)$ is prime we have the distribution relation
$$
\mathrm{Norm}_{K(\mathfrak{ml})/K(\mathfrak{m})}
\Theta(1;\mathfrak{ml},\fa)^e=\left\{\begin{array}{ll}
    \Theta(1;\mathfrak{m},\fa)& \mathrm{if\ }
    \mathfrak{l}\mid\mathfrak{m}\\
    \Theta(1;\mathfrak{m},\fa)^{1-\sigma_{\mathfrak{l}}^{-1}}
    &\mathrm{else}
\end{array}\right.
$$
where $e=w_{\mathfrak{m}}/w_{\mathfrak{ml}}$.  In particular, the
sequence $\Theta(1;\ff p^k,\fa)$ is norm compatible for $k>0$. We
denote by $\beta(\fa)$ and $\beta^*(\fa)$ the images of this sequence
in $U_\p(K(\ff p^\infty))$ and $U_{\p^*}(K(\ff p^\infty))$,
respectively.

Define $J\subset \Lambda(K(\ff p^\infty))_{\RR_0}$ to be the
annihilator of $\mu_{p^\infty}$.  Then $J$ is the ideal generated by
$\sigma_\mathfrak{b}-\mathrm{N}(\mathfrak{b})$ as $\mathfrak{b}$
ranges over integral ideals prime to $\ff p$.

\begin{Thm}\label{yager}(Yager)
  There are isomorphisms of $\Lambda(K(\ff
  p^\infty))_{\RR_0}$-modules
  $$
  U_\p(K(\ff p^\infty))_{\RR_0} \iso J, \quad U_{\p^*}(K(\ff
  p^\infty))_{\RR_0} \iso J,
  $$
  which take
  $$
  \beta(\fa)\mapsto
  (\sigma_\mathfrak{a}-\mathrm{N}(\mathfrak{a}))\cdot\mu_\p,
  \quad \beta^*(\fa)\mapsto
  (\sigma_{\fa}-\mathrm{N}(\fa))\cdot\mu_{\p^*}.
  $$
\end{Thm}
\begin{proof}
  This is Proposition III.1.4 of \cite{deShalit}.
\end{proof}


\subsection{The twisted Euler system} 
\label{euler system section}

 
We continue to denote by $\fa$ a nontrivial $\co_K$-ideal prime to $\ff p$.
If $\mathfrak{m}\in I(\fa)$, define a unit
$$
\vartheta_\fa(\mathfrak{m})= \mathrm{Norm}_{K(\mathfrak{m}\ff p)/
  K(\mathfrak{m})} \Theta(1;\mathfrak{m}\ff p,\fa).
$$
If we let $\mathcal{K}_\fa=\cup_{\mathfrak{m}\in I(\fa)}
K(\mathfrak{m})$, the elements
$$
\vartheta_\fa(\mathfrak{m})\in H^1(K(\mathfrak{m}),\Z_p(1)),
$$
with $\fa$ fixed, form an Euler system for $(\Z_p(1),\ff
p,\mathcal{K}_\fa)$ in the sense of \cite{rubin}.

\begin{Prop}\label{euler system}
  Let $\q=\p$ or $\p^*$.  There is an Euler system $c=c_{\fa}$ for
  $(T_\p,\ff p,\mathcal{K}_\fa)$ and an injection (the Coleman map)
  $$
  \mathrm{Col}_\q:\HH^1(K(\ff p^\infty)_\q,T_\p)_{\RR_0} \map{}
  \Lambda(K(\ff p^\infty))_{\RR_0}
  $$
  with the following property: if we set
  $$
  z=\mil c(K(\ff p^k)) \in \HH^1(K(\ff p^\infty), T_\p),
  $$
  and let $\loc_\q(z)$ be its image in $\HH^1(K(\ff
  p^\infty)_\q, T_\p)$, then $\mathrm{Col}_\q$ sends
  $$
  \loc_\q (z)\otimes 1\mapsto \mu_{\q}(K(\ff
  p^\infty),\psi_{\p^*},\fa).
  $$
  The image of $\mathrm{Col}_\q$ is $\Tw_{\psi_{\p^*}}(J)$, the
  ideal generated by all elements of the form
  $\sigma_\mathfrak{b}-\psi_\p(\sigma_\mathfrak{b})$ with
  $\mathfrak{b}$ prime to $\ff p$.
\end{Prop}
\begin{proof}
  This follows from the ``twisting'' theorems of Chapter 6 of
  \cite{rubin}.  The $G_K$-module $T_\p$ is isomorphic to the twist of
  $\Z_p(1)$ by the character $\omega^{-1}\psi_\p=\psi_{\p^*}^{-1}$,
  where $\omega$ is the cyclotomic character.  A choice of such an
  isomorphism determines an isomorphism of $\Z_p$-modules
  $$
  \HH^1(K(\ff p^\infty),\Z_p(1))\map{\phi}\HH^1(K(\ff
  p^\infty),T_\p)$$
  satisfying
  $\phi\circ\Tw_{\psi_{\p^*}}^{-1}(\lambda)=\lambda\circ\phi$ for any
  $\lambda\in \Lambda(K(\ff p^\infty))$.  Similarly, if $\q=\p$ or
  $\p^*$, there is an isomorphism of $\Z_p$-modules
  $$U_\q(K(\ff p^\infty))\iso\HH^1(K(\ff
  p^\infty)_\q,\Z_p(1))\map{\phi_\q} \HH^1(K(\ff p^\infty)_\q,T_\p)$$
  which is compatible with $\phi$ and the localization map $\loc_\q$.
  
  The Euler system $c_\fa$ is the twist of $\vartheta_\fa$ by
  $\psi_{\p^*}$, and by construction $$z=\phi(\mil\vartheta_\fa(\ff
  p^k))= \phi(\mil\Theta(1;\ff p^k,\fa)).$$
  See Section 6.3 of
  \cite{rubin}, especially Theorem 6.3.5.  In particular,
  $\loc_\q(z)=\phi_\q(\beta_\fa)$, and if we define $\mathrm{Col}_\q$
  to be the composition
  $$\HH^1(K(\ff p^\infty)_\q,T_\p)_{\RR_0}\map{\phi_\q^{-1}}
  U_\q(K(\ff p^\infty))_{\RR_0}\iso J\map{\Tw_{\psi_{\p^*}}}
  \Tw_{\psi_{\p^*}}(J)$$
  then $\mathrm{Col}_\q$ is an isomorphism of
  $\Lambda(K(\ff p^\infty))_{\RR_0}$-modules with the desired
  properties.
\end{proof}

\begin{Def} We say that the prime $p$ is \textit{anomalous} if $p$
divides 
$$
(1-\psi(\p))(1-\psi(\p^*)),
$$
or equivalently if there is any $p$-torsion defined over $\Z/p\Z$ on
the reduction of $E$ at $p$. 
\end{Def}

\begin{Lem}\label{local norms}
  Let $F\subset K_\infty$ contain $K$ and let $\q=\p$ or $\p^*$.  The
  natural corestriction map
\begin{equation}
\label{descent sequence}
\HH^1(K_{\infty,\q},T_\p)\otimes_{\Lambda(K_\infty)}\Lambda(F)
\map{}\HH^1(F_{\q},T_\p)
\end{equation} 
is an isomorphism if either $\q=\p^*$ or if $\q=\p$ and $p$ is not
anomalous.  If $F=C_m D_n$ for some $0\le m,n\le\infty$, then the map
is injective with finite cokernel.
\end{Lem}
\begin{proof}
  Let $v$ be any place of $K_\infty$ above $\p^*$, and denote also by
  $v$ the place of $F$ below it. By Lemma \ref{no local torsion},
  $E(K_{\infty,v})[\p^*]=0$, and the
  inflation-restriction sequence shows that
  $$
  H^1(F_v,W_{\p^*})\map{}
  H^1(K_{\infty,v},W_{\p^*})^{\Gal(K_{\infty,v}/F_v)} $$
  is an
  isomorphism.  This implies that the restriction map
  $$
  H^1(F_{\p^*},W_{\p^*})\map{}
  H^1(K_{\infty,\p^*},W_{\p^*})^{\Gal(K_\infty/F)}
  $$
  is an isomorphism. By local duality, the map
  $$
  \HH^1(K_{\infty,\p^*},T_\p)\otimes_{\Lambda(K_\infty)}\Lambda(F)
  \map{}\HH^1(F_{\p^*},T_\p)
  $$
  is an isomorphism.  If $p$ is not anomalous then $E(L)[\p]=0$ for
  any $p$-power extension $L/K_\p$, and the same argument as above
  shows that
  $$
  \HH^1(K_{\infty,\p},T_\p)\otimes_{\Lambda(K_\infty)}\Lambda(F)
  \map{}\HH^1(F_{\p},T_\p)
  $$
  is an isomorphism.
  
  Now suppose $F=C_m D_n$, let $v$ be a place of $K_\infty$ above
  $\p$, and suppose also that $p$ is anomalous; this implies that all
  $\p^*$-power torsion of $E$ is defined over $K_{\infty,v}$, and in
  fact is defined over the unique unramified $\Z_p$-extension of
  $K_v$.  Set $L=C_\infty D_n$, and define
  $$
B=H^0(L_v, W_{\p^*})=H^0(K_v^\mathrm{unr}\cap L_v, W_{\p^*}).
$$
  Note that $B$ is finite, since $K_v^\mathrm{unr}\cap L_v$ is a
  finite extension of $\Q_p$.  We have the inflation-restriction
  sequence
  $$
  H^1(L_v/F_v,B)\map{}H^1(F_v,W_{\p^*})\map{\res} H^1(
  L_v,W_{\p^*})^{\Gal(L_v/F_v)} \map{}H^2(L_v/F_v,B).
  $$
  Since $\Gal(L_v/F_v)\iso\Z_p$ is of cohomological dimension one, the
  final term of this sequence is trivial, and so the map $\res$ is
  surjective. The first term of the sequence is isomorphic to
  $B/(\gamma-1)B$ for any generator $\gamma$ of $\Gal(L_v/F_v)$, and
  this is finite since $B$ is.
  
  Now consider the restriction map
  $$
H^1(L_v,W_{\p^*})\map{}H^1(K_{\infty,v}, W_{\p^*})^{
    \Gal(K_{\infty,v}/L_v)}.
$$
  Since $K_{\infty,v}/L_v$ is a
  $\Z_p$-extension, the restriction is surjective, exactly as above.
  Similarly, the kernel is isomorphic to $W_{\p^*}/(\sigma-1)W_{\p^*}$
  for any generator $\sigma$ of $\Gal(K_{\infty,v}/L_v)$.  Such a
  $\sigma$ acts on $W_{\p^*}$ through some scalar $\not=1$, and since
  $W_{\p^*}$ is divisible the kernel of restriction is trivial.  As
  above, local duality gives the stated results.
\end{proof}

\begin{Prop}\label{descent}
  Let $F\subset K_\infty$ be a (possibly infinite) extension of $K$,
  and let $\q=\p$ or $\p^*$.  Suppose that one of the following holds
\begin{enumerate}
\item \label{descent one} $\q=\p^*$,
\item \label{descent two} $\q=\p$ and $p$ is not anomalous,
\item $\q=\p$ and $F=D_m C_n$ for $0\le m,n\le\infty$.
\end{enumerate}
There is an injection of $\Lambda_{\RR}(F)$-modules
$$
\HH^1(F_{\q},T_\p)_ \RR \map{}\Lambda(F)_\RR
$$
taking $\loc_{\q}(c(F))\otimes 1$ to $\mu_{\q}(F;\psi_{\p^*},\fa)$.
The image of this map is the ideal generated by 
$\lambda(F; \psi_{\p^*},\fa)$ as $\fa$ varies,
and if $F=D_\infty$ the map is an isomorphism.
If $p$ does not divide $[K(\ff):K]$ and if either (\ref{descent one})
or (\ref{descent two}) holds, the result is true with $\RR$ replaced by
$\RR_0$.
\end{Prop}
\begin{proof}
  Let $w$ be a place of $K(\ff p^\infty)$ above $p$, set $H=\Gal(K(\ff
  p^\infty)_w/K_{\infty,w} )$ and consider the inflation-restriction
  sequence
  $$
  H^1(H,W_{\p^*})\map{} H^1(K_{\infty,w},W_{\p^*})\map{}H^1(K(\ff
  p^\infty)_w, W_{\p^*})\map{} H^2(H,W_{\p^*}).
  $$
  The first and last terms are finite.  If $p$ does not divide
  $[K(\ff):K]$ then $p$ does not divide $|H|$, and so the kernel and
  cokernel of restriction are trivial.  Applying local duality and
  considering the semi-local cohomology, we see that the map
  $$
  \HH^1(K(\ff p^\infty)_\q,T_\p)\otimes_{\Lambda(K(\ff p^\infty))}
  \Lambda(K_\infty)\map{}\HH^1(K_{\infty,\q}, T_\p)
  $$
  has finite kernel and cokernel, and is an isomorphism if
  $[K(\ff):K]$ is prime to $p$.
  
  By Lemma \ref{local norms} the natural map
  $$
  \HH^1(K(\ff p^\infty)_\q,T_\p)\otimes_{\Lambda(K(\ff
    p^\infty))}\Lambda(F) \map{}\HH^1(F_{\q},T_\p)
  $$
  has finite kernel and cokernel, and so becomes an isomorphism
  upon applying first $\hat{\otimes}\RR_0$ and then $\otimes\RR$.
  Tensoring the Coleman map of Proposition \ref{euler system} with
  $\Lambda(F)$ yields the desired map.  If $F=D_\infty$
  then $F$ is disjoint from the extension of $K$ cut out by
  $\psi_\p$, 
  and it follows that the image of the map above is 
  the ideal of $\Lambda(F)_\RR$ generated by all elements of the 
  form $\gamma-\alpha$, where $\gamma$ runs over $\Gal(D_\infty/K)$
  and $\alpha$ runs over $\Z_p^\times$.  This ideal is all of 
  $\Lambda(F)_\RR$.
\end{proof}


\subsection{Main conjectures}
\label{main conjectures}


Throughout this subsection we fix a topological generator 
$\gamma\in\Gal(K_\infty/D_\infty)$,
and we let $I=(\gamma-1)\Lambda(K_\infty)$.  

Denote by $W$ the sign in the functional equation of $L(\psi,s)$.  Let
$\fa$ be an integral ideal of $\co_K$ prime to $p\ff$, and recall that
$\mathcal{K}_\fa$ is the union of all ray class fields of $K$ of
conductor prime to $\fa$.  Let $c_{\fa}$ be the Euler system for
$(T_\p, \ff p, \mathcal{K}_\fa)$ of Proposition \ref{euler system},
and for any $F\subset K_\infty$, let
$$
c_\fa(F)=\mil c_\fa(F')\in \HH^1(F, T_\p).
$$
be the limit as $F'$ ranges over subfields of $F$ finite over $K$.
Let $\CC_\fa(F)$ be the $\Lambda(F)$-submodule of
$\CS_\rel(F,T_\p)$ generated by $c_{\fa}(F)$, and let
$\CC(F)$ be the submodule generated by $\CC_\fa(F)$ as $\fa$
varies over all ideals prime to $p\ff$.
Define
$$
Z(F)=\CS_\rel(F,T_\p)/\CC(F),
$$
and define $Z_\fa(F)$ similarly, replacing $\CC$ by $\CC_\fa$.

\begin{Rem}
  It is clear from the definitions that
  $\Sel_\rel(F,T_\p)=H^1(F,T_\p)$ for every extension $F/K$.  In
  particular $c_\fa(F)\in \CS_\rel(F,T_\p)$.
\end{Rem}

\begin{Lem}\label{unramification}
  For every $F\subset\mathcal{K}_\fa$ finite over $K$, the class
  $c_\fa(F)$ is unramified at every prime of $F$ not dividing $p$.
\end{Lem}
\begin{proof}
  This follows from Corollary B.3.5 of \cite{rubin}, and the fact that
  the class $c_\fa(F)$ is a universal norm in the cyclotomic
  direction.
\end{proof}

\begin{Prop}\label{parity II}
  The submodule $\loc_{\p^*}(\CC(D_\infty)) \subset
  \HH^1(D_{\infty,\p^*},T_{\p})$ is nontrivial if and only if $W=1$.
  The submodule $\loc_{\p}(\CC(D_\infty)) \subset
  \HH^1(D_{\infty,\p},T_{\p})$ is nontrivial if and only if $W=-1$.
  In particular, $\CC(D_\infty)\not=0$ regardless of the value of $W$.
\end{Prop}
\begin{proof}
  Suppose $W=1$.  By Corollary \ref{parity},
  $$
  \mu_\p(D_\infty,\psi_{\p^*})=0\hspace{1cm}
  \mu_{\p^*}(D_\infty,\psi_{\p^*})\not=0.
  $$
  We may choose $\fa$ so that
  $\lambda(D_\infty;\psi_{\p^*},\fa)\not=0$, and then Proposition
  \ref{descent} implies that $\loc_{\p^*}(\CC_\fa(D_\infty))\not=0$.
  Therefore $\loc_{\p^*}(\CC(D_\infty))\not=0$.  
  On the other hand, for
  every choice of $\fa$, $\mu_{\p}(D_\infty;\psi_{\p^*},\fa)=0$, and
  so Proposition \ref{descent} implies that 
  $\loc_\p(\CC_\fa(D_\infty))$
  is trivial in
  $\HH^1(D_{\infty,\p},T_\p)_\RR$
  By Proposition \ref{local ranks},
  $\HH^1(D_{\infty,\p},T_\p)$ is torsion-free, and so
  $\loc_\p(\CC_\fa(D_\infty))$ is trivial in
  $\HH^1(D_{\infty,\p},T_\p)$.  The case $W=-1$ is entirely similar.
\end{proof}

Now armed with a nontrivial Euler system, we may apply the general
theory introduced by Kolyvagin and developed by 
Kato, Perrin-Riou, and Rubin.

\begin{Prop}\label{main conjecture one}
  The $\Lambda(D_\infty)$-module $X_\str^*(D_\infty)$ is torsion
  and $\CS_\rel(D_\infty,T_\p)$ is torsion free of rank one.
Furthermore, we have the divisibility of characteristic ideals
$$
\mathrm{char}(X_\str^*(D_\infty))\mathrm{\ \ divides\ \ }
\mathrm{char}(Z(D_\infty)).
$$
If $W=1$ then the $\Lambda(D_\infty)$-module $X^*(D_\infty)$ is torsion.
If $W=-1$ then $X^*(D_\infty)$ has rank one and 
$\CS(D_\infty,T_\p)=\CS_\rel(D_\infty,T_\p)$.

\end{Prop}
\begin{proof}
  By Proposition \ref{parity II}, we may choose some $\fa$ so that
  $\CC_\fa(D_\infty)\not=0$.  
The first claim follows from Theorem 2.3.2
  of \cite{rubin} (together with Lemma \ref{unramification} and the
  remarks of Section 9.2 of \cite{rubin}), and the second then follows
  from Lemmas \ref{equal ranks} and
  \ref{change of group} and Theorem \ref{strong selmer rank}.
Applying Theorem 2.3.3 of \cite{rubin} as $\fa$
varies over all integral ideals prime to $p\ff$, one 
obtains the divisibility of characteristic ideals.

Suppose $W=1$. The image of $\CC(D_\infty)$ under
  $$
  \loc_{\p^*}:\CS_\rel(D_\infty,T_\p)\map{}\HH^1(D_{\infty,\p^*},T_\p)
  $$
  is nontrivial by Proposition \ref{parity II}, and since both modules
  are torsion-free of rank one this map must be injective with 
torsion cokernel.  By
  Proposition \ref{global duality} we obtain the exact sequence
  $$
  0\map{}\CS_\rel(D_\infty,T_\p)\map{}\HH^1(D_{\infty,\p^*},T_\p)
  \map{}X^*(D_\infty)\map{}X^*_\str(D_\infty)\map{}0
$$
which shows that $X^*(D_\infty)$ is torsion.

Suppose $W=-1$.  By Proposition \ref{parity II} and the exact sequence
(\ref{tate one}), $\CC(D_\infty)\subset \CS(D_\infty,T_\p)$.  Since
$\CC(D_\infty)\not=0$ and
$\CS(D_\infty,T_\p)\subset\CS_\rel(D_\infty,T_\p)$, it follows that
$\CS(D_\infty,T_\p)$ is torsion free of rank one. From Lemma
\ref{change of group} and Proposition \ref{strong selmer rank}, we see that
$X^*(D_\infty)$ has rank one. Furthermore, the image of $\loc_{\p^*}$
in the exact sequence (\ref{tate one}) (with $F=D_\infty)$ must be a
torsion module, and hence must be trivial.  Therefore
$\CS(D_\infty,T_\p)=\CS_\rel(D_\infty,T_\p)$.
\end{proof}

Let $\mathcal{K}=K(E[p^\infty])$ and
define abelian extensions of $\mathcal{K}$ as follows:
\begin{itemize}
\item $\mathcal{M}_\rel$ is the maximal abelian pro-$p$-extension
of $\mathcal{K}$ unramified outside of $p$,
\item $\mathcal{M}$ is the maximal abelian pro-$p$-extension
of $\mathcal{K}$ unramified outside of $\p^*$,
\item $\mathcal{M}_\str$ is the maximal abelian pro-$p$-extension
of $\mathcal{K}$ unramified everywhere.
\end{itemize}
Let
$\mathcal{E}$ be the inverse limit of the groups $\co_F^\times\otimes\Z_p$ 
over subfields $F\subset \mathcal{K}$ containing $K$.
For any integral ideal $\fa$ of $K$ prime to $p\ff$, let
$\mathcal{U}_\fa\subset \mathcal{E}$ be the submodule generated by
the (untwisted) elliptic unit Euler system $\vartheta_\fa$ of
\S \ref{euler system section}.  Let $\mathcal{U}$ be the submodule
generated by all such $\mathcal{U}_\fa$. 
By Kummer theory we may view $\mathcal{U}\subset\mathcal{E}\subset
\HH^1(\mathcal{K},\Z_p(1))$.

The following result is essentially due to Coates (see \cite{coates}
Theorem 12).

\begin{Lem}\label{twisted selmer}
There is a group isomorphism
$$
\alpha:H^1(\mathcal{K},\Q_p/\Z_p)\iso H^1(\mathcal{K},W_{\p^*})
$$
satisfying $\alpha\circ \Tw_{\psi^*}(\lambda)=\lambda\circ \alpha$
for every $\lambda\in\Lambda(\mathcal{K})$, where 
$\Tw_{\psi_{\p^*}}:\Lambda(\mathcal{K})\map{}\Lambda(\mathcal{K})$
is the ring automorphism of \S \ref{L functions}.
This map restricts to an isomorphism (of groups, not 
$\Lambda(\mathcal{K})$-modules)
$$
\Hom(\Gal(\mathcal{M}/\mathcal{K}),\Q_p/\Z_p)\iso 
\Sel(\mathcal{K},W_{\p^*})
$$
and similarly for the relaxed and restricted Selmer groups, replacing
$\mathcal{M}$ by $\mathcal{M}_\rel$ and $\mathcal{M}_\str$, respectively.
Similarly, there is a group isomorphism
$$
\beta: \HH^1(\mathcal{K},\Z_p(1))\iso \HH^1(\mathcal{K},T_\p)
$$
satisfying $\beta\circ \Tw^{-1}_{\psi_{\p^*}}(\lambda)=\lambda\circ \beta$
for every $\lambda\in\Lambda(\mathcal{K})$.
This isomorphism identifies $\mathcal{E}$ with 
$\CS_\rel(\mathcal{K}, T_\p)$ and $\mathcal{U}$ with 
$\mathcal{C}(\mathcal{K})$.
\end{Lem}

\begin{proof}
The existence of 
$$
\alpha:H^1(\mathcal{K},\Q_p/\Z_p)\iso H^1(\mathcal{K},W_{\p^*})
$$
follows from the twisting theorems of \cite{rubin} \S 6.2, once one fixes 
an isomorphism
$W_{\p^*}\iso(\Q_p/\Z_p)(\psi_{p^*})$.  From the definitions,
together with Proposition \ref{global duality},
we have the following characterizations of our Selmer groups
in $H^1(\mathcal{K},W_{\p^*})$:
\begin{itemize}
\item
$\Sel_\rel(\mathcal{K},W_{\p^*})$ consists of the classes locally
trivial away from $p$,
\item
$\Sel(\mathcal{K},W_{\p^*})$ consists of the classes locally
trivial away from $\p^*$,
\item $\Sel_\str(\mathcal{K},W_{\p^*})$ consists of the classes 
everywhere locally trivial.
\end{itemize}
The isomorphism $\alpha$ identifies each of these Selmer groups
with the subgroup of classes in $H^1(\mathcal{K},\Q_p/\Z_p)$
satisfying the same local conditions, and so it suffices to 
check that for every place $v$ of $\mathcal{K}$, the condition 
``locally trivial at $v$'' agrees with the condition 
``unramified at $v$''.  For $v$ not dividing $p$
this is \cite{rubin} Lemma B.3.3, and the case $v|p$ is identical:
fix a place $v$ of $\mathcal{K}$ and
note that regardless of the rational prime below $v$,
$\mathcal{K}_v$ always contains the 
unique unramified $\Z_p$-extension of $K_v$.
In particular, if $\mathcal{K}^\mathrm{unr}$ is the maximal unramified
extension of $\mathcal{K}$, then $\Gal(\mathcal{K}^\mathrm{unr}/\mathcal{K})$
has trivial pro-$p$-part, and so
$H^1(\mathcal{K}^\mathrm{unr}/\mathcal{K},\Q_p/\Z_p)=0.$

For the compact cohomology group $\HH^1(\mathcal{K},T_\p)$ the
existence of $\beta$ is proved in the same fashion, using the fact
that $\psi_\p\psi_{\p^*}$ is the cyclotomic character.
That $\beta$ identifies the relaxed Selmer group with the unit
group is a consequence of the discussion above, since  local duality
shows that the unramified conditions agree with the relaxed conditions
everywhere locally.
The identification of $\mathcal{U}_\fa$ with $\CC_\fa(\mathcal{K})$ 
is immediate from the construction of the twisted Euler system
$c_\fa$ from the elliptic units $\vartheta_\fa$
in Proposition \ref{euler system}.
\end{proof}

Decompose
$$\Gal(\mathcal{K}/K)\iso \Delta\times\Gal(K_\infty/K)$$
where $\Delta\iso \Gal(K(E[p])/K)$, and let
$\psi_\q=\chi_\q\eta_\q$ be the associated decomposition of $\psi_\q$, for
$\q=\p$ of $\p^*$.

\begin{Cor}
We have the
equality of characteristic ideals in $\Lambda(K_\infty)$
$$
\Tw_{\eta_\p^*}^{-1}\big(
\mathrm{char}(\Gal(\mathcal{M}/\mathcal{K})^{\chi_{\p^*}})\big)
=
\mathrm{char}(X^*(K_\infty))
$$
and similarly for the relaxed and restricted Selmer groups.
Also,
$$
\Tw_{\eta_\p^*}^{-1}\big(
\mathrm{char}(\mathcal{E}/\mathcal{U})^{\chi_{\p^*}}\big)
=
\mathrm{char}(Z(K_\infty)).
$$
\end{Cor}
\begin{proof}
This follows easily by taking $\Delta$-invariants of the
$\Lambda(\mathcal{K}/K)$-modules of the Lemma above
(cf. e.g. \cite{rubin} Lemma 6.1.2).  One must remember our convention,
Remark \ref{twisted action}, about the $\Lambda(\mathcal{K})$-action
on $X^*(\mathcal{K})$.
\end{proof}

\begin{Thm}(Rubin)\label{rubin}
The $\Lambda(K_\infty)$-module $X^*(K_\infty)$ is torsion,
$\CS_\rel(K_\infty,T_\p)$ has rank one,
and 
\begin{eqnarray*}
\mathrm{char}(X^*(K_\infty))&=& 
\mathrm{char}\big(\HH^1(K_{\infty,\p^*},T_\p)/\loc_{\p^*}\CC(K_\infty)\big)
\\
\mathrm{char}(X^*_\str(K_\infty))&=&\mathrm{char}(Z(K_\infty)).
\end{eqnarray*}
\end{Thm}
\begin{proof}
In view of Lemma \ref{twisted selmer} and its corollary, this is a
twisted form of the main results of \cite{rubin91}.
\end{proof}

\begin{Rem}
The fact that $X^*(K_\infty)$ is torsion is originally due to Coates
\cite{coates}.
\end{Rem}

The following proposition follows from a deep result of Greenberg.
Strictly speaking, it is not needed to prove the main result of this
section, Theorem \ref{main conjecture two} below, but it is helpful for 
understanding the case $W=1$.  See Remark \ref{main comment}.

\begin{Prop}\label{null}
The $\Lambda(K_\infty)$-modules $X^*_\rel(K_\infty)$ and
$X^*(K_\infty)$ have no nonzero pseudo-null submodules.
\end{Prop}
\begin{proof}
It is a theorem of Greenberg \cite{greenberg0} that 
$\Gal(\mathcal{M}_\rel/\mathcal{K})$ has no non-zero pseudo-null submodules,
and so Lemma \ref{twisted selmer} implies that $X^*_\rel(K_\infty)$
also has none.  By \cite{PR84} \S II.2 Th\'eor\`eme 23, 
$X^*(K_\infty)$ also has no nontrivial pseudo-null submodules.
\end{proof}

\begin{Lem}\label{local control one}
For $\q=\p$ or $\p^*$, the kernel of the restriction map
$$
H^1(D_{\infty,\q},W_{\p^*})\map{}
H^1(K_{\infty,\q},W_{\p^*})
$$
is finite.
\end{Lem}
\begin{proof}
The kernel of the restriction map is isomorphic to
$$
H^0(K_{\infty,\q},W_{\p^*})/IH^0(K_{\infty,\q},W_{\p^*})
$$
by the inflation-restriction sequence and \cite{rubin}
Lemma B.2.8.
The finiteness follows from Lemma \ref{local norms} and local
duality.
\end{proof}

\begin{Lem}\label{local control two}
The semi-local restriction map
$$
\bigoplus_{w|\ff}H^1(D_{\infty,w},W_{\p^*})
\map{}
\bigoplus_{w|\ff}H^1(K_{\infty,w},W_{\p^*}).
$$
is injective.
\end{Lem}
\begin{proof}
This is \cite{PR84} II.7, Lemme 13 (or \cite{deShalit}, proof of  
Lemma IV.3.5),
together with the isomorphism
$$
H^1(L,W_{\p^*})\iso H^1(L,E)[\p^{*\infty}]
$$
for any algebraic extension $L/K_w$ (it suffices to prove this 
isomorphism for
finite extensions, where it is a consequence of the Kummer sequence
and the fact that $E(L)$ has a  finite
index pro-$\ell$ subgroup, where $\ell\not=p$ is the residue
characteristic of $w$). 
\end{proof}

We will need the following slight generalization of the control theorems
of Mazur and Perrin-Riou.

\begin{Prop}\label{global control}
The dual to the restriction map 
\begin{equation}\label{control restriction}
\Sel_\rel(D_\infty,W_{\p^*})\map{}
\Sel_\rel(K_\infty,W_{\p^*})[I]
\end{equation}
is an isomorphism of $\Lambda(D_\infty)$-modules
$$
X^*_\rel(K_\infty)/I X^*_\rel(K_\infty)\map{}X^*_\rel(D_\infty).
$$
The analogous maps for $X^*$ and $X^*_\str$ are surjective with
finite cokernel.
\end{Prop}
\begin{proof}
Let $S$ be the set of places of $K$ consisting of the archimedean place
and the primes dividing $p\ff$, and denote by $K_S/K$ the maximal
extension of $K$ unramified outside $S$.
By Lemma \ref{no local torsion} $H^0(K_S/K_\infty,W_{\p^*})=0$, and so the
inflation-restriction sequence shows that the restriction map
$$
H^1(K_S/D_\infty,W_{\p^*})\map{}H^1(K_S/K_\infty,W_{\p^*})[I]
$$
is an isomorphism.  
As in the proof of Lemma \ref{ramification}, for $w$ a prime of 
$D_\infty$ not lying above a prime of $S$
the local condition $H^1_f(D_{\infty,w},W_{\p^*})$ is equal to the 
unramified condition, and similarly for $K_\infty$. 
From the definition of the relaxed Selmer group, 
we have the commutative diagram with exact rows
$$
\xymatrix{
0\ar[r]  &  \Sel_\rel(D_\infty,W_{\p^*})\ar[r]\ar[d]  &  
H^1(K_S/D_\infty,W_{\p^*})\ar[r]\ar[d]
& \bigoplus H^1(D_{\infty,w},W_{\p^*})\ar[d] \\
0\ar[r]  &  \Sel_\rel (K_\infty,W_{\p^*})[I]\ar[r]  
&  H^1(K_S/K_\infty,W_{\p^*})[I]\ar[r]
& \bigoplus H^1(K_{\infty,w},W_{\p^*})
}
$$ where the direct sums are over places $w|\ff$.  In particular,
since the middle vertical arrow is an isomorphism, the restriction map
(\ref{control restriction}) is injective, and to bound the cokernel of
this map it suffices to bound the kernel of the right vertical arrow
in the diagram above.  This kernel is trivial by Lemma \ref{local
control two}.  This completes the proof for the relaxed Selmer groups.

In order to prove the result for the true Selmer groups, we replace the
top and bottom rows of the commutative diagram above with the exact
sequence (\ref{tate four}) applied with $F=D_\infty$ and $F=K_\infty$,
respectively.  Again by the snake lemma, it then suffices to bound the
kernel of restriction
$$
H^1(D_{\infty,\p},W_{\p^*})\map{}H^1(K_{\infty,\p},W_{\p^*}),
$$ 
and this is the content of Lemma \ref{local control one}.
Similarly, one deduces the result for the strict Selmer group from the
result for the true Selmer group by using the exact sequence
(\ref{tate two}), together with another application of Lemma
\ref{local control one}.
\end{proof}

\begin{Def}
We define the \emph{descent defect} $\defect\subset \Lambda(D_\infty)$ by
$$
\defect=\mathrm{char}_{\Lambda(D_\infty)}(X^*_\str(K_\infty)[I]).
$$
\end{Def}

\begin{Cor}\label{no zero}
The descent defect $\defect$ is nonzero, and
we have the equality of ideals in $\Lambda(D_\infty)$ 
$$
\mathrm{char}_{\Lambda(D_\infty)}(X^*_\str(D_\infty))=
\mathrm{char}_{\Lambda(K_\infty)}(X^*_\str(K_\infty))\cdot 
\defect
$$
If $W=1$ then 
$$
\mathrm{char}_{\Lambda(D_\infty)}(X^*(D_\infty))=
\mathrm{char}_{\Lambda(K_\infty)}(X^*(K_\infty)).
$$
\end{Cor}
\begin{proof}
By Proposition \ref{main conjecture one} $X^*_\str(D_\infty)$
is a torsion $\Lambda(D_\infty)$-module, and so by Proposition
\ref{global control}, the same is true of
$X_\str^*(K_\infty)/I X_\str^*(K_\infty)$.
The claim now follows from
\cite{rubin91} Lemma 6.2 (i).

When $W=1$, $X^*(D_\infty)$ is a torsion module and the proof is
identical, except that now \cite{rubin91} Lemma 6.2 (i) and Proposition
\ref{null} above show that $X^*(K_\infty)[I]=0$.
\end{proof}

Recall from Propositions \ref{parity II} and \ref{main conjecture one}
that $Z(D_\infty)$ is a torsion $\Lambda(D_\infty)$-module.
The following proposition gives the other half of the descent
from $K_\infty$ to $D_\infty$.

\begin{Prop}\label{compact control}
The natural maps of $\Lambda(D_\infty)$-modules 
\begin{eqnarray*}
\CS_\rel(K_\infty,T_\p)/I\CS_\rel(K_\infty,T_\p)
&\map{}&\CS_\rel(D_\infty,T_\p)\\
Z(K_\infty)/I Z(K_\infty)&\map{}&Z(D_\infty)
\end{eqnarray*}
are injective, and their cokernels have characteristic ideal $\defect$.
The same holds with $Z$ replaced by $Z_\fa$ for
any ideal $\fa$ prime to $p\ff$.
\end{Prop}
\begin{proof}
Let $L\subset K_\infty$ be finite over $K$. As always, let
$S$ the set of places of $K$ consisting of the archimedean place
and the primes dividing $p\ff$.  Let $K_S/K$ be the maximal
extension of $K$ unramified outside $S$. 
From the Poitou-Tate nine-term exact sequence we extract the
exact sequence
$$
0\map{}X^*_\str(L)\map{}H^2(K_S/L,T_\p)\map{}
\bigoplus_{v\in S}H^0(L_v,W_{\p^*})^\vee,
$$
(for example, by taking $B_v=0$ in \cite{perrin-riou} Proposition 4.1).
Passing to the limit as $L$ varies and taking $I$-torsion gives
$$
0\map{}X^*_\str(K_\infty)[I]\map{}
H^2(K_S/K,T_\p\otimes\Lambda(K_\infty))[I]\map{}
\bigoplus_{v|p\ff}H^0(K_{\infty,v},W_{\p^*})^\vee[I],
$$
where we have used Shapiro's lemma to identify
$$
H^2(K_S/K,T_\p\otimes\Lambda(K_\infty))\iso\mil H^2(K_S/L,T_\p).
$$
The final term in the exact sequence is the Pontryagin dual of 
\begin{eqnarray*}\lefteqn{
\bigoplus_{v|p\ff}H^0(K_{\infty,v},W_{\p^*})
/(\gamma-1)\bigoplus_{v|p\ff}H^0(K_{\infty,v},W_{\p^*}) }\hspace{2cm} \\
&\iso& \bigoplus_{v|p\ff} 
H^1(K_{\infty,v}/D_{\infty,v},H^0(K_{\infty,v},W_{\p^*}) ),
\end{eqnarray*}
which is the kernel of restriction 
$$
\bigoplus_{v|p\ff}H^1(D_{\infty,v},W_{\p^*})
\map{}\bigoplus_{v|p\ff}H^1(K_{\infty,v},W_{\p^*}).
$$
Lemmas \ref{local control one} and \ref{local control two}
(for $v|p$ and $v|\ff$, respectively) show that this kernel is finite.
By Corollary \ref{no zero} we conclude 
that $H^2(K_S/K,T_\p\otimes\Lambda(K_\infty))[I]$ is a torsion 
$\Lambda(D_\infty)$-module with characteristic ideal equal to $\defect$.

From the $\Gal(K_S/K)$-cohomology of 
$$
0\map{}T_\p\otimes\Lambda(K_\infty)\map{\gamma-1}
T_\p\otimes\Lambda(K_\infty)\map{}T_\p\otimes\Lambda(D_\infty)\map{}0
$$
we deduce that the map
$$
H^1(K_S/K,T_\p\otimes\Lambda(K_\infty))\otimes_{\Lambda(K_\infty)}
\Lambda(D_\infty) \map{}
H^1(K_S/K,T_\p\otimes\Lambda(D_\infty))
$$
is injective with torsion cokernel of characteristic ideal $\defect$.
Again using Shapiro's lemma,
together with Lemma \ref{ramification},
we see that the map
$$ 
\CS_\rel(K_\infty,T_\p)/I\CS_\rel(K_\infty,T_\p)\map{}\CS_\rel(D_\infty,T_\p)
$$
is injective with cokernel of characteristic ideal $\defect$. 
The map 
$$
\CC(K_\infty)/I\CC(K_\infty)\map{}
\CC(D_\infty)
$$
is visibly surjective, since this merely asserts that the twisted elliptic
units are universal norms in the cyclotomic direction.
The snake lemma now proves the claim.
\end{proof}

\begin{Prop}\label{strong main conjecture}
We have the equality of characteristic ideals
$$
\mathrm{char}(X_\str^*(D_\infty))=\mathrm{char}(Z(D_\infty)).
$$
\end{Prop}
\begin{proof}
Let $\fa$ be an ideal of $K$ prime to $p\ff$.
Using the fact that $I$ is 
principal, the snake lemma gives the exactness of
$$ 
\CS_\rel(K_\infty,T_\p)[I]\map{}Z_\fa(K_\infty)[I]
\map{}\CC_\fa(K_\infty)/I\CC_\fa(K_\infty).$$  
The leftmost term is trivial by Lemma \ref{equal ranks},
and the term on the right is isomorphic to $\Lambda(D_\infty)$,
since $\CC_\fa(K_\infty)$ is free of rank one over $\Lambda(K_\infty)$.
Therefore $Z_\fa(K_\infty)[I]$ is a torsion-free $\Lambda(D_\infty)$-module.
On the other hand, the quotient $Z_\fa(K_\infty)/IZ_\fa(K_\infty)$
is a torsion $\Lambda(D_\infty)$-module (by Propositions 
\ref{main conjecture one} and \ref{compact control}),
and so \cite{rubin91} Lemma 6.2 (i)
tells  us that  $Z(K_\infty)[I]$ is a torsion 
$\Lambda(D_\infty)$-module.  We conclude that $Z_\fa(K_\infty)[I]=0$.
Now by \cite{rubin91} Lemma 6.2 (ii),
$$
\mathrm{char}_{\Lambda(K_\infty)}(Z_\fa(K_\infty))\cdot \Lambda(D_\infty)
=\mathrm{char}_{\Lambda(D_\infty)}(Z_\fa(K_\infty)/IZ_\fa(K_\infty)).
$$
Applying Proposition \ref{compact control} gives
$$
\mathrm{char}_{\Lambda(K_\infty)}(Z_\fa(K_\infty))\cdot \defect
=\mathrm{char}_{\Lambda(D_\infty)}(Z_\fa(D_\infty)).
$$
Now let $\fa$ vary and apply
Theorem \ref{rubin} and Proposition \ref{no zero} to get
\begin{eqnarray*}
\mathrm{char}_{\Lambda(D_\infty)}(X_\str^*(D_\infty))&=&
\mathrm{char}_{\Lambda(K_\infty)}(X_\str^*(K_\infty))\cdot \defect\\
&=&\mathrm{char}_{\Lambda(K_\infty)}(Z(K_\infty))\cdot \defect
\end{eqnarray*}
proving the claim.
\end{proof}

\begin{Thm}\label{main conjecture two}
\

\begin{enumerate}
\item If $W=1$, then 
\begin{enumerate}
\item $\CS(D_\infty,T_\p)=0$,
\item $X^*(D_\infty)$ is a torsion $\Lambda(D_\infty)$-module, 
\item the ideal of $\Lambda_\RR(D_\infty)$ generated by 
$\mathrm{char}(X^*(D_\infty))$ is equal to the ideal
generated by the $p$-adic $L$-function
$\mu_{\p^*}(D_\infty,\psi_{\p^*})$.  If $p$ does not divide 
$[K(\ff):K]$ the same holds with $\RR$ replaced by $\RR_0$.
\end{enumerate}

\item  If $W=-1$, then
\begin{enumerate}
\item $\CS(D_\infty,T_\p)$ is a torsion-free
  $\Lambda(D_\infty)$-module of rank one,
\item $X^*(D_\infty)$ has rank one,
\item $\mathrm{char}\big(X^*_\str(D_\infty)\big)
  =\mathrm{char}\big(\CS(D_\infty,T_\p)/\CC(D_\infty)\big)$.
\end{enumerate}
\end{enumerate}
\end{Thm}
\begin{proof}
The first two claims of (1) and (2) all follow from 
Lemma \ref{equal ranks}, Theorem \ref{strong selmer rank}, and
Proposition \ref{main conjecture one}.
When $W=1$ the determination of the characteristic ideal 
follows from Propositions \ref{descent} and
\ref{strong main  conjecture}, using the exact sequence
$$
0\map{}Z(D_\infty)\map{}\HH^1(D_{\infty,\p^*},T_\p)/\loc_{\p^*}\CC(D_\infty)
\map{}X^*(D_\infty)\map{}X^*_\str(D_\infty)\map{}0.
$$
When $W=-1$ the claim follows from the final statement of
Proposition \ref{main conjecture one} and from Proposition
\ref{strong main conjecture}.
\end{proof}

\begin{Rem}\label{main comment}
The case $W=1$ can be deduced more directly from the first equality of
Theorem \ref{rubin} and the second part of Corollary \ref{no zero}
(which requires Proposition \ref{null}), by using Lemma \ref{local
norms}, the local analogue of Proposition \ref{compact control}.  This
avoids the application of the Euler system machinery directly over
$D_\infty$ (Proposition \ref{main conjecture one}) needed to prove
Proposition \ref{compact control}, or, more precisely, to prove the
nontriviality of $\defect$.  When $W=-1$, the ideals appearing in the
first equality of Theorem \ref{rubin} have trivial image in
$\Lambda(D_\infty)$, so one seems to have no recourse but to prove
some form of Proposition \ref{compact control}.
\end{Rem}


\section{The $p$-adic height pairing}
\label{height section}


Throughout this section we assume $W=-1$ and we set
$$
\Delta=\Gal(D_\infty/K)\hspace{1cm}\Gamma=\Gal(C_\infty/K).
$$
We will frequently identify $\Gamma\iso\Gal(K_\infty/D_\infty)$.  For
any nonnegative integer $n$, let $\Delta_n=\Delta/\Delta^{p^n}$ and
similarly for $\Gamma$.  Let $\mathcal{I}$ be the kernel of the
natural projection $\Lambda(C_\infty)\map{}\Lambda(K)$ and set
$\aug=\mathcal{I}/ \mathcal{I}^2$.  Many authors use some choice of
``logarithm'' $\lambda:\Gamma\map{}\Z_p$ to define the $p$-adic height
pairing.  Following the fashion of the day, we instead take
$$
\lambda:\Gamma\map{}\aug
$$
to be the isomorphism $\gamma\mapsto
\gamma-1$, and so obtain a $\aug$-valued height pairing.

\subsection{The linear term}

Choose a generator $\gamma\in\Gamma$
and fix some integral ideal
$\fa\subset\co_K$ prime to $\ff p$.  For every $K\subset L\subset
K_\infty$ we set
$$c_\fa(L)=\mil c_\fa(L'),$$
where $c_\fa$ is the Euler system of
Proposition \ref{euler system} and the limit is taken over all
subfields $L'\subset L$ finite over $K$.  We may identify
$$\Lambda(K_\infty)_{\RR_0}\iso\Lambda(D_\infty)_{\RR_0}[[\Gamma]]$$
and expand $\mu_{\p^*}(K_\infty,\psi_{\p^*},\fa)$ as a power series in
$\gamma-1$,
$$\mu_{\p^*}(K_\infty,\psi_{\p^*},\fa)=\mathcal{L}_{\fa,0}+
\mathcal{L}_{\fa,1}(\gamma-1)+ \mathcal{L}_{\fa,2}(\gamma-1)^2+\cdots.
$$
Similarly we may expand
$$\mu_{\p^*}(K_\infty,\psi_{\p^*})=\mathcal{L}_{0}+
\mathcal{L}_{1}(\gamma-1)+ \mathcal{L}_{2}(\gamma-1)^2+\cdots.
$$
By Corollary \ref{parity} and the assumption that $W=-1$, we have
$$\mathcal{L}_{0}=\mu_{\p^*}(D_\infty,\psi_{\p^*})=0.$$
It follows
that also $\mathcal{L}_{\fa,0}=0$. By Proposition \ref{descent} the
image of $c_\fa(D_\infty)$ in $\HH^1(D_{\infty,\p^*},T_\p)$ is
trivial, and so by Proposition \ref{global duality},
$c_\fa(D_\infty)\in \CS(D_\infty,T_\p)$.

\begin{Lem}\label{alg derivative}
  Set $F_n=D_nC_\infty$.  For every $n$ there is a unique element
$$
\beta_{n}\in \HH^1(F_{n,\p^*},T_\p)_\RR
$$
such that
$$
(\gamma-1)\beta_{n}=\loc_{\p^*}\big( c_\fa(F_n)\big).
$$
Let $\alpha_{n}$ be the image of $\beta_{n}$ in
$\HH^1(D_{n,\p^*},T_{\p})_\RR$. The elements $\alpha_n$ are
norm-compatible and they define an element $\alpha_{\infty} \in
\HH^1(D_{\infty,\p^*},T_\p)_{\RR}$.  The Coleman map of Proposition
\ref{descent} identifies
$$
\HH^1(D_{\infty,\p^*},T_\p)_{\RR}\iso\Lambda(D_\infty)_{\RR}
$$
and takes $\alpha_\infty$ to $\mathcal{L}_{\fa,1}$.
\end{Lem}
\begin{proof}
This is immediate from Proposition \ref{descent}, the fact that
$\mathcal{L}_{\fa,0}=0$, and the definition of $\mathcal{L}_{\fa,1}$.
\end{proof}

Tate local duality defines a pairing
$$
\langle\ ,\ \rangle_n:\HH^1(D_{n,\p^*},T_{\p})\times
\HH^1(D_{n,\p^*},T_{\p^*})\map{}\Z_p
$$
whose kernel on either side is
the $\Z_p$-torsion submodule, and the induced pairing on the quotients
by the torsion submodules is perfect.

The height pairing of the following theorem has been studied by many
authors, including Mazur-Tate, Nekov\'{a}\v{r}, Perrin-Riou, and
Schneider. The fourth property of the pairing, the height formula, is
due to Rubin, and plays a crucial role in what follows.

\begin{Thm}\label{height one}
  For every nonnegative integer $n$ there is a canonical (up to sign)
  $p$-adic height pairing
  $$h_n : \Sel(D_n,T_\p)\times \Sel(D_n,T_{\p^*})\map{}\Q_p\otimes\aug$$
  satisfying the following properties
\begin{enumerate}
\item there is a positive integer $k$, independent of $n$, such that
  $h_n$ takes values in $p^{-k}\Z_p\otimes\aug$
\item if $a\in \Sel(D_n,T_\p)$, $b\in \Sel(D_n,T_{\p^*})$, and
  $\sigma\in \Delta_n$, then
  $$h_n(a^\sigma, b^\sigma)=h_n(a,b)$$
\item if $a_n\in \Sel(D_n,T_\p)$, $b_{n+1}\in \Sel(D_{n+1},T_{\p^*})$,
  and $\res$ and $\cor$ are the restriction and corestriction maps
  relative to $D_{n+1}/D_n$, then
  $$h_{n+1}(\res(a_n), b_{n+1})=h_n(a_n,\cor(b_{n+1}))$$
\item (height formula) for every $b \in \Sel(D_n,T_{\p^*})$, we have
  (up to sign)
  $$h_n\big( c_\fa(D_n), b \big)=
  \langle\alpha_{n},\loc_{\p^*}(b)\rangle_n\otimes (\gamma-1).$$
\end{enumerate}
\end{Thm}
\begin{proof}
This will be proved in the next section.
\end{proof}

Define the $\Lambda(D_\infty)$-adic Tate pairing
$$\langle\ ,\ \rangle_\infty:\HH^1(D_{\infty,\p^*},T_\p)_\RR
\otimes_{\Lambda(D_\infty)_\RR}
\HH^1(D_{\infty,\p^*},T_{\p^*})^\iota_\RR\iso\Lambda(D_\infty)_\RR$$
by
$$
\langle a_\infty, b_\infty\rangle_\infty= \mil
\sum_{\sigma\in\Delta_n}\langle a_n^\sigma,b_n\rangle_n\cdot
\sigma^{-1}
$$
and define the $\Lambda(D_\infty)$-adic height pairing
$$
h_\infty : \CS(D_\infty,T_\p)_\RR\otimes_{\Lambda(D_\infty)_\RR}
\CS(D_\infty,T_{\p^*})^\iota_\RR\map{}
\Lambda(D_\infty)_\RR\otimes_{\Z_p}\aug
$$
similarly. The element $\alpha_{\infty}\in
\HH^1(D_{\infty,\p^*},T_\p)$ satisfies
\begin{equation}\label{rubin height}
h_\infty(c_\fa(D_\infty),b_\infty)=\langle \alpha_{\infty},
\loc_{\p^*}(b_\infty)\rangle_\infty\otimes(\gamma-1)
\end{equation}
for every $b_\infty\in \CS(D_\infty,T_{\p^*})$.

\begin{Def} \label{D:reg}
  Define the \emph{anticyclotomic regulator}, $\reg$, to be the 
characteristic ideal of the cokernel of $h_\infty$.  
\end{Def}

Define $\reg(\CC_\fa)$ to be the characteristic ideal of 
the cokernel of 
$$h_\infty|_{\CC_\fa}:\CC_\fa\otimes_{\Lambda(D_\infty)_\RR}
\CS(D_\infty,T_{\p^*})^\iota\map{}\Lambda(D_\infty)_\RR\otimes_{\Z_p}\aug,
$$ 
and let $\eta$ be  the ideal
  $$\eta=\mathrm{char}\big(\HH^1(D_{\infty,\p^*},
  T_{\p^*})/\loc_{\p^*}(\CS(D_\infty,T_{\p^*}))\big).$$
  From Proposition \ref{main conjecture one} and the results of
  Section \ref{anticyclotomic theory} we have that
  $\CS_\str(D_\infty,T_{\p^*})$ is trivial.
  The exactness of (\ref{tate three}) then shows that $\eta\not=0$. 

\begin{Prop}\label{height}
  There is an equality of ideals in $\Lambda(D_\infty)_\RR$
\begin{eqnarray*}
\reg\cdot\mathrm{char}\big(\CS(D_\infty,T_\p)/\CC_{\fa}\big)&=&
\reg(\CC_\fa)\\ &=&(\mathcal{L}_{\fa,1})\cdot \eta^\iota.
\end{eqnarray*}
\end{Prop}
\begin{proof}
The first equality is clear.
The height formula (\ref{rubin height}) implies that
the image of $h_\infty|_{\CC_\fa}$ is equal to
$$
\langle\alpha_\infty,\loc_{\p^*}
(\CS(D_\infty,T_{\p^*})_\RR)^\iota\rangle_\infty
\otimes\aug\subset \Lambda(D_\infty)_\RR\otimes_{\Z_p}\aug
$$
The second equality now follows from Lemma 
\ref{alg derivative} and the fact that
the $\Lambda(D_\infty)$-adic Tate pairing is an isomorphism.
\end{proof}

\begin{Thm} \label{T:bsd}
  Let $\mathcal{X}$ denote the ideal of $\Lambda(D_\infty)_\RR$
  generated by the characteristic ideal of the
  $\Lambda(D_\infty)$-torsion submodule of $X^*(D_\infty)$.  Then we
  have the equality of ideals
  $$\reg\cdot \mathcal{X}= (\mathcal{L}_1).$$
\end{Thm}
\begin{proof}
  If we replace $\p$ by
  $\p^*$ and take $F=D_\infty$ in the second pair of exact sequences
  of Proposition \ref{global duality}, we obtain the exact sequence
  $$0\map{}\HH^1(D_{\infty,\p^*},
  T_{\p^*})/\loc_{\p^*}(\CS(D_\infty,T_{\p^*}))
  \map{}X_\rel(D_\infty)\map{}
  X(D_\infty)\map{}0.$$
  Taking $\Lambda(D_\infty)$-torsion and
  applying Lemma \ref{change of group} and Theorem \ref{strong
  selmer rank} we obtain
  $$\mathrm{char}\big(X_\str^*(D_\infty)\big)=\eta^\iota\cdot
  \mathcal{X}.
  $$
  Letting $\fa$ vary in Proposition \ref{height}, the claim follows from
  Theorem \ref{main conjecture two} (2), part (c).
\end{proof}

 
\subsection{The height formula}


In this section we sketch Perrin-Riou's construction of the $p$-adic
height pairing
$$
h_n : \Sel(D_n,T_\p)\times \Sel(D_n,T_{\p^*})\map{}\Q_p\otimes\aug
$$
of Theorem \ref{height one}, as well Rubin's proof of the height
formula. Our exposition closely follows that of \cite{rubin94}, to
which we refer the reader for details.

For every $0\le k\le\infty$, let $L_k=D_n C_k$.

\begin{Lem}
  Fix a place $v$ of $D_n$ and some extension of it to $L_\infty$. The
  submodule of $H^1_\f(D_{n,v},T_\p)$ defined by
  $$
H^1_\f(D_{n,v},T_\p)^\mathrm{univ}= \bigcap_{k}\cor\ 
  H^1_\f(L_{k,v},T_\p)
$$
  has finite index, and the index is bounded as
  $v$ and $n$ vary.
\end{Lem}
\begin{proof}

First assume that $v$ divides $\p^*$, so that $H^1_\f(D_{n,v},T_\p)$
is the torsion submodule of $H^1(D_{n,v},T_\p)$, which is in turn
isomorphic to $H^0(D_{n,v},W_\p)$.  It clearly suffices to show that 
this is bounded as $n$ varies.  But $K_v(W_\p)$ is unramfied, and 
$D_{\infty,v}$ is a ramified $\Z_p$-extension of $K_v$, so this is 
clear.

From now on suppose that $v$ does not divide $\p^*$, so that
$$H^1_\f(L_{k,v},T_\p)=H^1(L_{k,v},T_\p).$$
By local duality, it suffices to bound the kernel of restriction
$$H^1(D_{n,v},W_{\p^*})\map{}H^1(L_{\infty,v},W_{\p^*}),$$
which is $H^1(L_{\infty,v}/D_{n,v}, M)$, where 
$M=E(L_{\infty,v})[\p^{*\infty}].$

If $v$ does not divide $\p$ then  $L_{\infty,v}$
is the unique unramified $\Z_p$-extension of $K_v$ (in particular 
it does not vary with $n$).
If $E$ does not have any $\p^*$-torsion defined over $K_v$,
then $K_v(E[\p^*])/K_v$ is a nontrivial extension of degree
dividing $p-1$, so $E$ has no $\p^*$-torsion defined over any 
$p$-extension of $K_v$,  and so $M=0$.
Assume $E[\p^*]$ is defined over $K_v$, and that $E$ has good reduction
at $v$.  Then $K_v(E[\p^{*\infty}])=L_{\infty,v}$ and so 
$M=W_{\p^*}$ is $p$-divisible.  If $\gamma$ is a generator of 
$\Gal(L_{\infty,v}/D_{n,v})$ then 
$$H^1(L_{\infty,v}/D_{n,v}, M)\iso M/(\gamma-1)M.$$
Since $\gamma$ acts as a scalar $\not=1$ on $M$, this group is trivial.
If $E$ has bad reduction at $v$ then $M$ is finite 
by the criterion of N\'eron-Ogg-Shafarevich, and the order of $M$
does not vary with $n$ (since $L_{\infty,v}$ does not vary).

Now assume that $v$ divides $\p$. The extension of $K_v$ 
generated by $E[\p^{*\infty}]$ is unramified, and since 
$K_v^\mathrm{unr}\cap L_{\infty,v}$ is a finite extension of $K_v$,
$M$ is finite.  If $\gamma$ generates $\Gal(L_{\infty,v}/D_{n,v})$
then using the exactness of
$$
0\map{}M^{\gamma=1}\map{}M\map{\gamma-1}M\map{}M/(\gamma-1)M\map{}0
$$
we see that the order of $H^1(L_{\infty,v}/D_{n,v},M)$
is equal to the order of 
$$
M^{\gamma=1}=E(D_{n,v})[\p^{*\infty}].
$$
Since $D_{\infty,v}$ is a ramified $\Z_p$-extension of $K_v$, this order
is bounded as $n$ varies.
\end{proof}

Fix $a\in \Sel(D_n,T_\p)$ and $b\in \Sel(D_n,T_{\p^*}).$ In view of
the above lemma, it suffices to define the height pairing of $a$ and
$b$ under the assumption that both are everywhere locally contained in
$H^1_\f(D_{n,v},T_\p)^\mathrm{univ}$.  Viewing $b$ as an element of
the larger group $H^1(D_n,T_{\p^*})$, $b$ defines an extension of
Galois modules
$$
0\map{}T_{\p^*}\map{}M_b^*\map{}\Z_p\map{}0,
$$
and taking
$\Z_p(1)$-duals we obtain an exact sequence
\begin{equation}\label{extension}
0\map{}\Z_p(1)\map{}M_b\map{}T_\p\map{}0.
\end{equation}
If $L/D_n$ is any finite extension, we may consider the global and
local Galois cohomology
$$
\xymatrix{ H^1(L,\Z_p(1)) \ar[r]\ar[d] & H^1(L,M_b)
  \ar[r]^{\pi_L}\ar[d]
  & H^1(L,T_\p) \ar[r]^{\delta_L}\ar[d] & H^2(L,\Z_p(1))\ar[d]\\
  H^1(L_w,\Z_p(1)) \ar[r] & H^1(L_w,M_b) \ar[r]^{\pi_{L_w}} &
  H^1(L_w,T_\p) \ar[r]^{\delta_{L_w}} & H^2(L_w,\Z_p(1)).  }
$$

\begin{Lem}\label{extension lemma}
  Let $L$ be a finite Galois extension of $D_n$ and suppose $a'\in
  H^1(L,T_{\p})$ satisfies $\cor(a')=a$. Then $a'$ is in the image of
  $\pi_L$.  For every place $w$ of $L$, $H^1_\f(L_w,T_\p)$ is
  contained in the image of $\pi_{L_w}$
\end{Lem}
\begin{proof}
  Let $\res$ be the restriction map from $D_n$ to $L$.  The connecting
  homomorphism $\delta_L$ is given (up to sign) by $\cup \res(b)$.  If
  $w$ is any place of $L$ and $v$ is the place of $D_n$ below it,
  $$\loc_w(\delta(a')) = \loc_w(a'\cup \res(b)) =\loc_v (a\cup b)=0,$$
  since $a$ and $b$ are everywhere locally orthogonal under the Tate
  pairing.  Thus $\delta_L(a')$ is everywhere locally trivial, and by a
  fundamental fact of class field theory it is globally trivial.  The
  proof of the second claim is similar.
\end{proof}

Class field theory gives a homomorphism
$$
\rho:\mathbf{A}_{D_n}^\times\map{}\Gal(D_nC_\infty/D_n)\iso
\Gamma\map{\lambda}\aug
$$
and we factor $\rho=\sum_v\rho_v$, the sum over all places of $D_n$.
By local Kummer theory we may view $\rho_v$ as a map
$$
\rho_v: H^1(D_{n,v},\Z_p(1))\map{}\aug,
$$
which can also be described as follows: the homomorphism
$$
\Gal(D_nC_\infty/D_n)\iso\Gamma\map{\lambda}\aug
$$
defines a class $\lambda \in H^1(D_n,\aug)$ (we always regard $\aug$ as
having trivial Galois action), and cup product with $\loc_v(\lambda)$
defines a map
$$
H^1(D_{n,v},\Z_p(1))\map{}H^2(D_{n,v},\aug(1))\iso\aug
$$
which
agrees (up to sign) with $\rho_v$.

Taking $L=D_n$ and $a'=a$ in  Lemma \ref{extension lemma}, we may choose some
$y^\glob\in H^1(D_n,M_b)$ with $\pi_{D_n}(y^\glob)=a$.  Fix a place
$v$ of $D_n$ and an extension of $v$ to $L_\infty$, and for every $k$
choose $y_{k,v}\in H^1_\f(L_{k,v},T_\p)$ which corestricts to
$\loc_v(a)$.  By Lemma \ref{extension lemma} we may choose some $y'_{k,v}\in
H^1(L_{k,v},M_b)$ such that $\pi_{L_{k,v}}(y'_{k,v})=y_{k,v}$.  Let
$\cor (y'_{k,v})$ be the image of $y'_{k,v}$ in $H^1(D_{n,v},M_b)$.
Then $\loc_v(y^\glob)-\cor(y'_{k,v})$ comes from some $w_{k,v}\in
H^1(D_{n,v},\Z_p(1))$, and we define
$$h_n(a,b)=\lim_{k\to\infty} \sum_v \rho_v(w_{k,v})$$
This limit
exists and is independent of all choices made.

We now sketch the proof of the height formula.  Suppose that
$a=c_\fa(D_n)$, and set $a_k=c_\fa(L_k)$,
$$a_\infty=\mil a_k\in \mil \Sel_\rel(L_k,T_\p).$$
By Lemma
\ref{extension lemma} there is a sequence $z_k\in H^1(L_k,M_b)$ with
$\pi_{L_k}(z_k)=a_k$.  Working semi-locally above $\q=\p$ or $\p^*$,
we have defined in the preceeding paragraph a sequence $y'_{k,\q}\in
H^1( L_{k,\q}, M_b)$ which lifts $y_{k,\q}$. The image of
$$t_{k,\q}=\loc_\q(z_k)- y'_{k,\q}\in H^1(L_{k,\q},M_b)$$
in
$H^1(D_{n,\q},M_b)$ comes from some $s_{k,\q}\in
H^1(D_{n,\q},\Z_p(1))$.

\begin{Prop}
  With notation as above
  $$h_n(a,b)=\lim_{k\to\infty}[\rho_\p (s_{k,\p})+
  \rho_{\p^*}(s_{k,\p^*})].$$

\end{Prop}
\begin{proof}
  This is Proposition 5.3 of \cite{rubin94}.
\end{proof}

Define $H_{k,\q}$ by the exactness of
$$0\map{}H_{k,\q}\map{}H^1(L_{k,\q},T_\p)\map{\cor}H^1(D_{n,\q},T_\p).$$
In \cite{rubin94} Section 4, one finds the definition of a derivative
operator
$$\Der_{k,\q}:H_{k,\q}\map{}H^1(D_{n,\q},T_\p/p^k T_\p) \otimes
\aug.$$
From the definition of $t_{k,\q}$, it is immediate that
$\pi_{L_{k,\q}}(t_{k,\q})\in H_{k,\q}$.  We set
$$t'_{k,\q}=\Der_{k,\q}(\pi_{L_{k,\q}}(t_{k,\q})) \in
H^1(D_{n,\q},T_\p/p^k T_\p) \otimes \aug.
$$
Proposition 4.3 of \cite{rubin94} then reads
$$
\lambda \cup s_{k,\q}=\delta_{D_{n,\q}}(t'_{k,\q})\in
H^2(D_{n,\q},(\Z/p^k\Z) (1))\otimes\aug$$
and so up to sign
\begin{eqnarray}\label{rubin cup}
\rho_\q(s_{k,\q})&\equiv& \mathrm{inv}_\q(\delta_{D_{n,\q}}(t'_{k,\q}) )
\pmod{p^k}\\
&=&\nonumber \mathrm{inv}_\q( t'_{k,\q}\cup\loc_\q(b))
\end{eqnarray}
where $\mathrm{inv}_\q$ is the semi-local invariant
$$H^2(D_{n,\q},(\Z/p^k\Z) (1))\otimes\aug\map{}\aug/p^k\aug.$$
From
the definition of $\Sel(D_n,T_{\p^*})$, we see that $\loc_\p(b)$ is a
torsion element.  If $p^\ell \loc_\p(b)=0$ then (\ref{rubin cup})
implies that $\rho_\p(s_{k,\p})$ is divisible by $p^{k-\ell}$.
Letting $k\to\infty$, we have $\rho_\p(s_{k,\p})\to 0$, leaving
\begin{equation}\label{shortening}
  h_n(a,b)=\lim_{k\to\infty}\rho_{\p^*}(s_{k,\p^*}).
\end{equation}

\begin{Lem}\label{limit formula}
  Suppose $d_k\in H_{k,\p^*}$ is such that $d_k\cup z=
  \loc_{\p^*}(a_k)\cup z$ for every $z\in
  H^1_\f(L_{k,\p^*},T_{\p^*})$.  Then for any sequence $x_k\in
  H^1_\f(L_{k,\p^*},T_{\p^*})$ such that the image of $x_k$ in
  $H^1_\f(D_{n,\p^*},T_{\p^*})$ is constant, 
  $$\lim_{k\to\infty}\mathrm{inv}_{\p^*}( \Der_{k,\p^*}(d_k)\cup x_0)=
  \lim_{k\to\infty}\sum_{\sigma\in\Gamma_k} \langle
  \loc_{\p^*}(a_k),\sigma\cdot x_k\rangle_{L_k,\p^*}\otimes \lambda (\sigma)
  $$
  where $\lambda$ is viewed as a character $\Gamma_k\map{}
  \aug/p^k\aug$.
\end{Lem}
\begin{proof}
  This is Lemma 5.1 of \cite{rubin94}.
\end{proof}

Recall that Lemma \ref{alg derivative} provides, for some choice of
generator $\gamma\in\Gamma$, a $\beta\in\HH^1(L_{\infty,\p^*},T_\p)$
such that $$(\gamma-1)\beta=\loc_{\p^*}(a_\infty).$$ Let $\alpha$ be
the image of $\beta$ in $H^1(D_{n,\p^*},T_{\p})$.  Write
$$\beta=\mil\beta_k\in\mil H^1(L_{k,\p^*},T_\p)$$
so that $\beta_0=\alpha$.  Fix some 
sequence $x_k\in H^1_\f(L_{k,\p^*},T_{\p^*})$ lifting
$\loc_{\p^*}(b)$.
Applying Lemma \ref{limit formula} with 
$d_k=\pi_{L_{k,\p^*}}(t_{k,\p^*})$ and
comparing with (\ref{rubin cup}) and (\ref{shortening}) gives
\begin{eqnarray*}
h_n(a,b)&=&\lim_{k\to\infty}\rho_{\p^*}(s_{k,\p^*})\\ 
\nonumber&=&\lim_{k\to\infty}\sum_{\sigma\in\Gamma_k}
\langle \loc_{\p^*}(a_k),x_k^\sigma\rangle_{L_k,\p^*}\otimes \lambda(\sigma)\\
&=&\lim_{k\to\infty}\sum_{\sigma\in\Gamma_k}
\langle (\gamma-1)\beta_k,\sigma\cdot x_k
\rangle_{L_k,\p^*}\otimes \lambda(\sigma)\\
&=& \lim_{k\to\infty}\sum_{i=1}^{p^k-1}
\langle (\gamma-1)\beta_k, i\gamma^i\cdot x_k\rangle_{L_k,\p^*}\otimes
(\gamma-1)\\
&=&\lim_{k\to\infty}\langle \beta_k, 
\mathrm{Norm}_{\Gamma_k}x_k\rangle_{L_k,\p^*}\otimes
(\gamma-1)\\
&=&\langle \alpha, \loc_{\p^*}(b)\rangle_{D_{n,\p^*}}
\otimes(\gamma-1).
\end{eqnarray*}
This completes the proof of the height formula.

\appendix
\numberwithin{Thm}{section}
\section{Proof of Theorem \ref{T:linterm}\\  by Karl Rubin}

In this appendix we prove Theorem \ref{T:linterm} of the introduction.
Essentially what we need to prove is that the anticyclotomic regulator
$\R$ of Definition \ref{D:reg} is nonzero.  The key tool is Theorem
\ref{bertrandthm} of Bertrand below, which says that on a CM elliptic
curve the $\p$-adic height of a point of infinite order is nonzero.
This is much weaker than saying that the $p$-adic height is
nondegenerate, but it suffices for our purposes.

We assume throughout this appendix that the sign in the functional equation 
of $L(E/\Q,s)$ is $-1$.  

We need to consider a slightly more general version of $p$-adic heights 
than appears in the main text.  If $F$ is a finite extension of $K$ in 
$K_\infty$, then there is a $p$-adic height pairing
$$
h_F : \Sel(F,T_\p) \otimes \Sel(F,T_{\p^*}) \to \Gal(K_\infty/F)\otimes\Qp.
$$
We are interested in three specializations of this pairing.  Namely, 
\begin{align*}
h_{F,\cycl} : \Sel(F,T_\p) \otimes \Sel(F,T_{\p^*}) &\to
\Gal(FC_\infty/F)\otimes\Qp, \\ 
h_{F,\anti} : \Sel(F,T_\p) \otimes \Sel(F,T_{\p^*}) &\to
\Gal(FD_\infty/F)\otimes\Qp, \\ 
h_{F,\p} : \Sel(F,T_\p) \otimes \Sel(F,T_{\p^*}) &\to
\Gal(FL_\infty/F)\otimes\Qp 
\end{align*}
where $L_\infty$ is the unique $\Zp$-extension of $K$ which is
unramified outside of $\p$, are defined by restricting the image of
$h_F$ to the appropriate group.  If $P \in E(F)$, we will write
$h_F(P) = h_F(a,b)$ where $a$ is the image of $P$ in $\Sel(F,T_\p)$
and $b$ is the image of $P$ in $\Sel(F,T_{\p^*})$.

The $p$-adic height pairing $h_n$ of Theorem \ref{height one} is the
composition of $h_{D_n,\cycl}$ with the isomorphism $\Gal(D_n
C_\infty/D_n) \cong \Gal(C_\infty/K) \cong \J$.

\begin{Thm}[Bertrand \cite{bertrand}]
\label{bertrandthm}
Suppose $F$ is a number field and $P \in E(F)$ is a point of infinite order.  
Then the $\p$-adic height $h_{F,\p}(P)$ is nonzero.
\end{Thm}

\begin{Lem}
\label{dependent}
Suppose $a \in \Sel(F,T_\p)$ and $b \in \Sel(F,T_{\p^*})$.  If two of 
$h_{F,\cycl}(a,b)$, $h_{F,\anti}(a,b)$, $h_{F,\p}(a,b)$ are zero, 
then so is the third.
\end{Lem}

\begin{proof}
Since $\Gal(K_\infty/F) \cong \Zp^2$, the projections 
$\Gal(K_\infty/F) \to \Gal(FC_\infty/F)$, 
$\Gal(K_\infty/F) \to \Gal(FD_\infty/F)$, and
$\Gal(K_\infty/F) \to \Gal(FL_\infty/F)$ 
are linearly dependent.  It follows that each of the three heights 
is a linear combination of the other two.
\end{proof}

\begin{Def}
For every $n$ define a submodule of $\Sel(D_n,T_\p)$, the 
{\em universal norms}, by
$$
\Sel(D_n,T_\p)^\univ = \bigcap_{m > n} \cor_{D_m/D_n} \Sel(D_m,T_\p).
$$
Define $\Sel(D_n,T_{\p^*})^\univ$
similarly.
\end{Def}

\begin{Lem}
\label{deg}
For every $n$ we have 
$$h_{D_n,\anti}(\Sel(D_n,T_\p)^\univ \otimes \Sel(D_n,T_{\p^*})^\univ) = 0.$$
\end{Lem}

\begin{proof}
This is a basic property of the $p$-adic height (see for example 
Proposition 4.5.2 of \cite{mazurtate}).
\end{proof}

\begin{Prop}
\label{controlprop}
The natural maps
$$
X(D_\infty) \otimes \Lambda(D_n) \to X(D_n)
$$
have kernel and cokernel which are finite and bounded independently of $n$.  
\end{Prop}

\begin{proof}
This is the standard ``Control Theorem'', see for example \cite{control}. 
\end{proof}

Recall that $\Delta_n = \Gal(D_n/K)$.

\begin{Lem}
\label{weak}
For every $n$, $\Sel(D_n,T_\p)^\univ \otimes \Qp$ and 
$\Sel(D_n,T_{\p^*})^\univ \otimes \Qp$ are 
free of rank one over $\Qp[\Delta_n]$.
\end{Lem}

\begin{proof}
This is proved in exactly the same way as Theorem 4.2 of
\cite{mrgrowth} (which proves a stronger statement about universal
norms in Mordell-Weil groups, assuming that all relevant
Tate-Shafarevich groups are finite), using Theorem \ref{main
conjecture two}. For completeness we give a proof here.

The exact sequence $0 \to E[\p^k] \to W_\p \to W_\p \to 0$ and the 
fact that $E(D_n) \cap E[\p] = 0$ show that 
$\Sel(D_n,E[\p^k]) = \Sel(D_n,W_\p)[\p^k]$ for every $k$.  
This and the Control Theorem (Proposition \ref{controlprop}) give us maps
\begin{multline*}
\Sel(D_n,T_\p) = \varprojlim_k \Sel(D_n,E[\p^k]) 
    = \varprojlim_k \Hom(\Sel(D_n,E[\p^k])^\wedge,\Z/p^k\Z) \\
    = \varprojlim_k \Hom(X(D_n),\Z/p^k\Z)
    = \Hom(X(D_n),\Zp) \\
    = \Hom_{\Lambda(D_n)}(X(D_n),\Lambda(D_n)) 
    \to \Hom_{\Lambda(D_\infty)}(X(D_\infty),\Lambda(D_n)) 
\end{multline*}
with finite kernel and cokernel bounded independently of $n$.  
This gives the bottom row, and passing to the inverse limit over $n$ gives 
the top row, of the commutative diagram with horizontal isomorphisms
$$
\begin{CD}
\CS(D_\infty,T_\p)\otimes\Qp @>\sim>> 
    \Hom_{\Lambda(D_\infty)}(X(D_\infty),\Lambda(D_\infty))\otimes\Qp \\
@VVV @VVV \\
\Sel(D_n,T_\p)\otimes\Qp @>\sim>> 
    \Hom_{\Lambda(D_\infty)}(X(D_\infty),\Lambda(D_n))\otimes\Qp
\end{CD}
$$
By Theorem \ref{main conjecture two}, the upper modules are free of
rank one over $\Lambda(D_\infty)\otimes\Qp$.  The kernel of the
right-hand vertical map is $I_n
\Hom_{\Lambda(D_\infty)}(X(D_\infty),\Lambda(D_\infty))\otimes\Qp$
where $I_n$ is the kernel of the map $\Lambda(D_\infty) \to
\Lambda(D_n)$, and so the kernel of the left-hand vertical map is $I_n
\CS(D_\infty,T_\p)\otimes\Qp$.  Hence the image of the left-hand
vertical map is free of rank one over $\Qp[\Delta_n]$.  But that image
is precisely $\Sel(D_n,T_\p)^\univ \otimes \Qp$.

The proof for $\Sel(D_n,T_{\p^*})^\univ \otimes \Qp$ is the same.
\end{proof}

\begin{Prop}
\label{zero}
If the anticyclotomic regulator $\R$ is zero, then the $p$-adic 
height pairing $h_n$ is identically zero on 
$\Sel(D_n,T_\p)^\univ \otimes \Sel(D_n,T_{\p^*})^\univ$.
\end{Prop}

\begin{proof}
Suppose $n \ge 0$, $\bfa = (a_n) \in \CS(D_\infty,T_\p)$ and $\bfb = (b_n)
\in \CS(D_\infty,T_\p^*)$.  Recall that $\Delta_n = \Gal(D_n/K)$.
Using property (3) of Theorem \ref{height one} and the definition of
$h_\infty$, we see that projecting $h_\infty(\bfa,\bfb)$ to $\Zp[\Delta_n]
\otimes \J$ gives $\sum_{\sigma\in\Delta_n}
h_n(a_n^\sigma,b_n)\sigma^{-1}$.

Now suppose $v_n \in \Sel(D_n,T_\p)^\univ$ and $v_n^* \in
\Sel(D_n,T_{\p^*})^\univ$.  Since $v_n$ and $v_n^*$ are universal
norms we can choose $\bfa = (a_n) \in \CS(D_\infty,T_\p)$ and $\bfb =
(b_n) \in \CS(D_\infty,T_\p^*)$ with $a_n = v_n$ and $b_n = v_n^*$.
If $\R = 0$, then $h_\infty(\bfa,\bfb) = 0$, and projecting to
$\Zp[\Delta_n] \otimes \J$ shows that $h_n(v_n,v_n^*) = 0$.
\end{proof}

\begin{Prop}
\label{universal}
If $n$ is sufficiently large then there are points $P \in E(D_n)$ of
infinite order such that the image of $P$ in $\Sel(D_n,T_\p)$ lies in
$\Sel(D_n,T_\p)^\univ$ and the image of $P$ in $\Sel(D_n,T_{\p^*})$
lies in $\Sel(D_n,T_{\p^*})^\univ$.  I.e., there are $D_n$-rational
points of infinite order which are universal norms in the Selmer
group.
\end{Prop}

\begin{proof}
Let $\psi$ denote the Hecke character of $K$ attached to $E$, so that 
$L(E/\Q,s) = L(\psi,s)$.  Fix an integer $n$.  

Choose a character $\chi$ of $\Delta_n$ of order $p^n$ (any two such
characters are conjugate under $G_\Q$, so the choice will not matter).
The Hecke $L$-function $L(\psi\chi,s)$ is the $L$-function of a
modular form $f_\chi$ on $\Gamma_0(Np^{2n})$, where $N$ is the
conductor of $E$.  Let $A_n$ denote the simple factor over $\Q$ of the
Jacobian $J_0(Np^{2n})$ corresponding to $f_\chi$.  Comparing
$L$-functions we see that there is an isogeny of abelian varieties
over $K$
$$
A_n \times \Res_{D_{n-1}/K} E  \sim \Res_{D_n/K} E
$$ 
where $\Res$ stands for the restriction of scalars.  Passing to
Mordell-Weil groups we get
$$
(A_n(K)\otimes \Q) \times (E(D_{n-1})\otimes \Q) 
\cong E(D_n)\otimes \Q, 
$$
and therefore, if $\gamma$ is a topological generator of $\Delta_\infty$,
\begin{equation}
\label{isog2}
A_n(K)\otimes\Q \cong (1-\gamma^{p^{n-1}})E(D_n)\otimes\Q.
\end{equation}

It follows from our assumption about the sign in the functional
equation of $L(E/\Q,s)$ that $L(\psi\rho,1) = 0$ for every character
$\rho$ of finite order of $\Delta_\infty$. By a theorem of Rohrlich
\cite{rohrlich}, there are only finitely many characters $\rho$ of
$\Delta_\infty$ such that the derivative $L'(\psi\rho,1) = 0$.
Suppose now that $n$ is large enough so that
\begin{itemize}
\item[(a)]
$L'(\psi\chi,1) \ne 0$, 
\item[(b)] 
$\mathrm{char}(X(D_\infty)_\tors)$ is relatively prime to 
$(\gamma^{p^n}-1)/(\gamma^{p^{n-1}}-1)$.  
\end{itemize}
Since $L'(\psi\chi,1) \ne 0$, the theorem of Gross and 
Zagier \cite{grosszagier} shows that 
\begin{equation}
\label{two}
\rank_{\Z} A_n(\Q) \ge \dim\,A_n = p^n - p^{n-1}.
\end{equation}
On the other hand, using \eqref{isog2} and the Control Theorem 
(Proposition \ref{controlprop}) we get
\begin{multline}
\label{three}
\rank_\Z A_n(\Q) \le \rank_{\Zp}(1-\gamma^{p^{n-1}})X(D_n) \\
    = \rank_{\Zp}(1-\gamma^{p^{n-1}})X(D_\infty)/(1-\gamma^{p^n})X(D_\infty)
\end{multline} 
Since $X(D_\infty)$ has $\Lambda(D_\infty)$-rank one (Theorem
\ref{main conjecture two}(2)), we conclude from condition (b) on $n$
that
$$
\rank_{\Zp}(1-\gamma^{p^{n-1}})X(D_\infty)/(1-\gamma^{p^n})X(D_\infty) 
    = p^n - p^{n-1}.
$$
It follows that we must have equality in \eqref{two} and \eqref{three}, and 
$$
\dim_{\Qp} (1-\gamma^{p^{n-1}})\Sel(D_n,T_\p) \otimes \Qp 
    =\rank_{\Zp}(1-\gamma^{p^{n-1}})X(D_n) = p^n - p^{n-1}.
$$
Since $(1-\gamma^{p^{n-1}})\Qp[\Delta_n]$ is a simple $\Qp[\Delta_n]$-module 
it follows that 
$$
(1-\gamma^{p^{n-1}}) \Sel(D_n,T_\p) \otimes \Qp \cong 
    (1-\gamma^{p^{n-1}})\Qp[\Delta_n]
$$
By Lemma \ref{weak} we now see that 
$$
(1-\gamma^{p^{n-1}}) \Sel(D_n,T_\p)^\univ \otimes \Qp
    = (1-\gamma^{p^{n-1}}) \Sel(D_n,T_\p) \otimes \Qp.
$$ 
In particular if $P$ is any point of infinite order in
$(1-\gamma^{p^{n-1}})E(D_n)$ (and we know that such points exist by
\eqref{isog2}) then some multiple of the image of $P$ in
$\Sel(D_n,T_\p)$ lies in $\Sel(D_n,T_\p)^\univ$.  In exactly the same
way some multiple of the image of $P$ in $\Sel(D_n,T_{\p^*})$ lies in
$\Sel(D_n,T_{\p^*})^\univ$, and the proposition is proved.
\end{proof}

\begin{proof}[Proof of Theorem \ref{T:linterm}]
Using Proposition \ref{universal}, find an $n$ and a point of 
infinite order $P \in E(D_n)$ whose images in $\Sel(D_n,T_\p)$ 
and $\Sel(D_n,T_{\p^*})$ are universal norms.  

By Bertrand's Theorem \ref{bertrandthm}, we have $h_{D_n,\p}(P) \ne 0$.  
By Lemma \ref{deg} we have $h_{D_n,\anti}(P) = 0$, and therefore by 
Lemma \ref{dependent} we have that $h_{D_n,\cycl}(P)$ (and hence $h_n(P)$) 
is nonzero.  

It now follows from Proposition \ref{zero} that the anticyclotomic
regulator $\R$ is nonzero, and so by Theorem \ref{T:bsd} the leading
term $\LL_1$ is nonzero. This is Theorem \ref{T:linterm}.
\end{proof}

\begin{Rem}
In the notation of Proposition \ref{universal}, the abelian variety
$A_n$ is isogenous to its twist by the quadratic character of $K/\Q$,
and so there are isomorphisms
$$
A_n(K)\otimes\Q\iso (\Res_{K/\Q}A_n)(\Q)\otimes\Q
\iso (A_n\times A_n)(\Q)\otimes\Q.
$$
Thus
\begin{eqnarray*}
p^n-p^{n-1} &=& 
\frac{1}{2}\ \rank_{\Z}A_n(K)\\
&=& \rank_{\co_K} E(D_n) -\rank_{\co_K} E(D_{n-1})
\end{eqnarray*}
for $n\gg 0$. This implies that there is a constant $c$ such that 
$$
\rank_{\co_K}E(D_n)= p^n+ c
$$ 
for $n\gg 0$.  The same asymptotic formula holds for the 
corank of the $\p$-primary Selmer group 
(by Theorem \ref{main conjecture two} and Proposition
\ref{controlprop}), and so the $\co_\p$-corank of 
of $\mbox{\cyr Sh}(E/D_n)_{\p^\infty}$ is bounded as $n$ varies.
\end{Rem}



\begin{thebibliography}{1}

\bibitem{arnold} T.  ~Arnold.
\newblock
Anticyclotomic main conjectures for CM modular forms.
\newblock Preprint. 2005.

\bibitem{bertrand} D.\ Bertrand. Propri\'et\'es arithm\'etiques
de fonctions th\^eta \`a plusieurs variables.  In: Number theory,
Noordwijkerhout 1983, {\em Lecture Notes in Math.}  {\bf 1068}
Springer: Berlin (1984) 17--22.

\bibitem{coates}
J. ~Coates.
\newblock Infinite descent on elliptic curves with complex multiplication,
in {\em Arithmetic and Geometry, Vol. I,} 
107--137, Progr. Math., 35, Birkhäuser. (1983).

\bibitem{deShalit}
E.~de~Shalit.
\newblock {\em {I}wasawa Theory of Elliptic Curves with Complex
  Multiplication}.
\newblock Academic Press, (1987).

\bibitem{greenberg0}
R. ~Greenberg.
\newblock On the structure of certain Galois groups.
\newblock {\em Invent. Math.} 47 (1978) 85--99.


\bibitem{greenberg}
R.~Greenberg.
\newblock On the {B}irch and {S}winnerton-{D}yer conjecture.
\newblock {\em Invent. Math.} 72 (1983) 241--265.

\bibitem{grosszagier}
B.\ Gross, D.\ Zagier, Heegner points and derivatives of $L$-series.
{\em Invent.\ Math.} {\bf 84} (1986) 225--320.


\bibitem{me}
B. ~Howard.
\newblock The Iwasawa theoretic Gross-Zagier theorem.  \emph{Comp. Math.} 141, No. 4 (2005), 811--846.

\bibitem{lang}
S.~Lang.
\newblock {\em Algebraic Number Theory}, Second Edition, Springer
Verlag, (1994).

\bibitem{martinet}
J. ~Martinet.
\newblock Character theory and Artin $L$-functions.
\newblock In: {\em Algebraic Number Fields}, A. Fr\"ohlich (ed.),
Academic Press, (1977) 1--88.

\bibitem{control}
   B.\ Mazur, Rational points of abelian varieties with values in
   towers of number fields.  {\em Invent.\ Math.} {\bf 18} (1972) 183--266.

\bibitem{mazur}
B. ~Mazur.
\newblock Modular curves and arithmetic.
\newblock In: {\em Proceedings of the International Congress of
Mathematicians (Warsaw, 1983)}, PWN, Warsaw (1984) 195--237.

\bibitem{mrgrowth}
   B.\ Mazur, K.\ Rubin, Studying the growth of Mordell-Weil.
   In: {\em Documenta math.}\ Extra Volume: Kazuya Kato's Fiftieth Birthday
   (2003) 585--607.

\bibitem{mazur-rubin}
B. ~Mazur and K. ~Rubin.
\newblock Kolyvagin Systems.
\newblock {\em Memoirs of the AMS,} {\bf 168}, (2004).

\bibitem{MRtwo}
B. ~Mazur and K. ~Rubin.
\newblock Elliptic curves and class field theory.
\newblock In: {\em  Proceedings of the International Congress of
Mathematicians (Beijing, 2002)}.

\bibitem{mazurtate}
   B.\ Mazur, J.\ Tate, Canonical height pairings via biextensions.
   In: Arithmetic and geometry, Vol.\ I, {\em Progr. Math.} {\bf 35}
   Birkhuser: Boston (1983) 195--237.

\bibitem{PR84}
B. ~Perrin-Riou.
\newblock Arithm\'etique des courbes elliptiques et th\'eorie
d'Iwasawa.
\newblock {\em Bull. Soc. Math. France, Memoires}, N.S. 17 (1984).


\bibitem{BPR1}
B. ~Perrin-Riou.
\newblock Fonctions $L$ $p$-adiques, th\'eorie d'Iwasawa et points de
Heegner.
\newblock {\em Bull. Soc. math. France}, 115 (1987) 399-456.

\bibitem{perrin-riou}
B.~Perrin-Riou.
\newblock Th\'{e}orie d'{I}wasawa et hauteurs $p$-adiques.
\newblock {\em Invent. Math.} 109 (1992) 137--185.

\bibitem{rohrlich}
D. ~Rohrlich.
\newblock On $L$-functions of elliptic curves and anticyclotomic
towers.
\newblock {\em Invent. Math.,} 75 (1984) 383--408.

\bibitem{rubin91}
K.~Rubin.
\newblock The ``main conjectures'' of 
{I}wasawa theory for imaginary quadratic
  fields.
\newblock {\em Invent. Math.} 103 (1991) 25--68.

\bibitem{rubin92}
K. ~Rubin.
\newblock $p$-adic {$L$}-functions and rational points on elliptic curves
with complex multiplication.
\newblock {\em Invent. Math.} 107 (1992) 323-350.

\bibitem{rubin94}
K.~Rubin.
\newblock Abelian varieties, $p$-adic heights and derivatives.
\newblock In {\em Algebra and Number Theory}. 
Walter de Gruyter and Co. (1994) 247--266.


\bibitem{rubin99}
K.~Rubin.
\newblock Elliptic curves with complex multiplication and the conjecture of
  Birch and Swinnerton-dyer.
\newblock In {\em Arithmetic Theory of Elliptic Curves 
(Cetraro, 1997)}, number
  1716 in Lecture Notes in Math., pages 167--234. Springer, (1999).

\bibitem{rubin}
K.~Rubin.
\newblock {\em Euler Systems}.
\newblock Princeton University Press, (2000).

\bibitem{weil}
A. ~Weil.
\newblock {\em Automorphic Forms and Dirichlet Series}, 
\newblock Lecture Notes in Math. No. 189, Springer, (1971).

\bibitem{yager}
R.~Yager.
\newblock On two variable $p$-adic {$L$}-functions.
\newblock {\em Ann. Math.} 115 (1982) 411--449.

\end{thebibliography}
\end{document}